# SUPERDIFFUSIVITY FOR A BROWNIAN POLYMER IN A CONTINUOUS GAUSSIAN ENVIRONMENT

By Sérgio Bezerra,[1] Samy Tindel and Frederi Viens[2]

*Université de Nancy, Université de Nancy and Purdue University*

This paper provides information about the asymptotic behavior of a one-dimensional Brownian polymer in random medium represented by a Gaussian field $W$ on $\mathbb{R}_+ \times \mathbb{R}$ which is white noise in time and function-valued in space. According to the behavior of the spatial covariance of $W$, we give a lower bound on the power growth (wandering exponent) of the polymer when the time parameter goes to infinity: the polymer is proved to be superdiffusive, with a wandering exponent exceeding any $\alpha < 3/5$.

**1. Introduction.** This paper is concerned with a model for a one-dimensional directed Brownian polymer in a Gaussian random environment (random medium) which can be briefly described as follows: the polymer itself, in the absence of any random environment, will simply be modeled by a Brownian motion $b = \{b_t; t \geq 0\}$, defined on a complete filtered probability space $(\mathcal{C}, \mathcal{F}, (\mathcal{F}_t)_{t \geq 0}, (P_b^x)_{x \in \mathbb{R}})$, where $P_b^x$ stands for the Wiener measure starting from the initial condition $x$. The corresponding expected value will be denoted[3] by $E_b^x$, or simply by $E_b$, when $x=0$. One may assume that $\mathcal{C}$ is the space of continuous functions started at 0.

The random environment will be represented by a centered Gaussian field $W$ indexed by $\mathbb{R}_+ \times \mathbb{R}$, defined on another complete probability space $(\Omega, \mathcal{G}, \mathbf{P})$ independent of $b$'s canonical space. Denoting by $\mathbf{E}$ the expected

Received March 2007; revised March 2007.
[1]Supported in part by CAPES.
[2]Supported in part by NSF Grant 0204999.

[3]This notation, which employs a subscript $b$ in a somewhat abusive way to indicate that an average with respect to the distribution of $b$ is taken, is now common in random medium theory, and has the advantage of reminding the reader that the randomness being averaged out is that of the Brownian $b$, not of the medium.

value with respect to **P**, the covariance structure of $W$ is given by

$$\mathbf{E}[W(t,x)W(s,y)] = [t \wedge s]Q(x-y), \tag{1.1}$$

for a given homogeneous covariance function $Q:\mathbb{R} \to \mathbb{R}$ satisfying some growth conditions that will be specified later on. In particular, the function $t \mapsto [Q(0)]^{-1/2}W(t,x)$ is a standard Brownian motion for any fixed $x \in \mathbb{R}$, and for every fixed $t \in \mathbb{R}_+$, the process $x \mapsto t^{-1/2}W(t,x)$ is a homogeneous Gaussian field on $\mathbb{R}$ with covariance function $Q$.

Once $b$ and $W$ are defined, the polymer measure itself can be described as follows: for any $t > 0$, the energy of a given path (or configuration) $b$ on $[0,t]$, under the influence of the random environment $W$, is given by the *Hamiltonian*

$$-H_t(b) = \int_0^t W(ds, b_s). \tag{1.2}$$

A completely rigorous meaning for this integral will be given in the next section, but for the moment, notice that for any fixed path $b$, $H_t(b)$ is a centered Gaussian random variable with variance $tQ(0)$. Based on this Hamiltonian, for any $x \in \mathbb{R}$, and a given constant $\beta$ (interpreted as the inverse of the temperature of the system), we define our (random) polymer measure $G_t^x$ (with $G_t := G_t^0$) as follows:

$$dG_t^x(b) = \frac{e^{-\beta H_t(b)}}{Z_t^x} dP_b^x(b) \quad \text{with } Z_t^x = E_b^x[e^{-\beta H_t(b)}]. \tag{1.3}$$

After early results in the Mathematical Physics literature (see [7] and [12]), links between martingale theory and directed polymers in random environments were established in [2] and [1], and over the last few years, several papers have shed some light on different types of polymer models: the case of random walks in discrete potential is treated, for instance, in [3], the case of Gaussian random walks is in [13] and [14], and the case of Brownian polymers in a Poisson potential is considered in [6]. On the other hand, the second author of this paper has undertaken in [16] the study of the polymer measure $G_t$ defined by (1.3). This latter model, which is believed to behave similarly to the other directed polymers mentioned above, has at least one advantage, from our point of view: it can be tackled with a wide variety of methods, some of which are new to the field: scaling invariances for both $b$ and $W$, stochastic analysis, Gaussian tools. Our long-term goal is to exploit such tools in order to get a rather complete description of the asymptotic behavior of the measure $G_t$.

In the present article we undertake this task by investigating the so-called *wandering exponent* $\alpha$, which measures the growth of the polymer when $t$ tends to $\infty$, and can be defined informally by the fact that, under the measure $G_t$, $\sup_{s \leq t}|b_s|$ behave like $t^\alpha$ for large times $t$. This kind of exponent has been studied in different contexts in [6, 13, 14, 15] and [18], yielding the conclusion that, for a wide number of models in dimension one, we should have $3/5 \leq \alpha \leq 3/4$. The true exponent conjectured by physicists is $\alpha = 2/3$.

Our understanding, from references [8, 10] and [11], is that physicists have come to this conjecture in dimension one, based on simulations (e.g., [10]) and on theoretical evidence as well as physical heuristics (in [8] where $\alpha$ is denoted by $\zeta$). The lower bound $\alpha \geq 3/5$ is confirmed mathematically in partially discrete settings (e.g., [14]). Our Section 3 provides an explanation of how our quantitative results confirm that $\alpha$ should be no less than $3/5$ if the environment's spatial memory, that is, its spatial correlation range, is short enough (cubic decay rate), and that superdiffusivity ($\alpha > 1/2$) is only guaranteed if this memory is not too long (decay rate exponent exceeding $5/2$). These long-spatial-memory situations are ones which do not seem to be considered in the mathematical or physical literature, so it is possible that the conjecture $\alpha = 2/3$ may not apply, although at this stage we have no evidence of any example of an upper bound result implying $\alpha < 2/3$.

In this paper we will see that, for our model, we have $\alpha \geq 3/5$. More specifically, we will prove the following.

THEOREM 1.1. *Let $\beta$ be any strictly positive real number. Assume that $Q : \mathbb{R} \to \mathbb{R}$ defined by (1.1) is a symmetric positive function, decreasing on $\mathbb{R}_+$ and such that, for some constant $\theta > 0$,*

$$(1.4) \qquad Q(x) = O\left(\frac{1}{|x|^{3+\theta}}\right) \qquad as \ x \to \pm\infty.$$

*In particular, $Q(0) < \infty$, which implies that $W$ defined in (1.1) is function-valued in $x$. Then, for any $\varepsilon > 0$, we have*

$$(1.5) \qquad \lim_{t \to \infty} \mathbf{P}\left[\frac{1}{t^{3/5-\varepsilon}} \left\langle \sup_{s \leq t} |b_s| \right\rangle_t \geq 1\right] = 1,$$

*where $\langle \cdot \rangle_t$ denotes expectation with respect to the polymer measure $dG_t^x(b)$ in (1.3).*

Our proof of this result inspires itself with some of the steps of Peterman's work in [14], where the same kind of growth bound has been established for a random walk in a Gaussian potential. Notice that, beyond generalizing his work from discrete to continuous space, we have been able to extend Petermann's result to a wider class of environments: indeed, we prove the relation (1.5) holds as soon as $Q$ satisfies the mild correlation decay assumption (1.4); Peterman assumed an exponential decay for $Q$. Moreover, many arguments had to be changed in order to pass from the random walk to the Brownian case. Having said all this, we must express our debt to Peterman's work which, unfortunately, has not been published beyond this Ph.D. dissertation [14] as directed by Erwin Bolthausen.

From the physical standpoint, it is worth noting that the above superdiffusivity theorem (wandering exponent $\alpha > 1/2$), which obviously does not hold for $\beta = 0$ (absence of random environment), holds nonetheless for all $\beta > 0$, that is, all temperatures. This is in contrast to the notion of strong disorder, defined and described at the end of the next section, a concept that we will study in detail in a separate publication. However, taken in a naive and intuitive sense, strong disorder is morally implied by superdiffusivity; lower bounds on wandering exponents that exceed $1/2$ thus appear as a convenient quantitative way of measuring this disorder, which is proved here to hold uniformly for all temperatures.

This paper is structured as follows. Section 2 defines the random environment $W$ and the Hamiltonian $H_t(b)$ rigorously, and discusses the relation between our wandering exponent $\alpha$ and the concept of strong disorder. Section 3 discusses the meaning of our main technical Hypothesis 1.4, what happens when one tries to weaken it, and a related open problem on the interplay between superdiffusivity and random environment correlation range. Section 3 also presents the main strategy for proving Theorem 1.1. The remainder of the paper is devoted to proving this theorem. Section 4 calculates the asymptotic correlation structure of space-time averages of $W$. Section 5 calculates similar asymptotics describing the interaction between $b$ and $W$. Section 6 presents an application of Girsanov's theorem for $b$ which estimates the penalization needed to force distant portions of $b$ back near the origin. Finally, with all these quantitative tools in hand, the proof of the theorem is completed in Section 7, which also contains a detailed heuristic description of this part of the proof.

The authors of this paper express their thanks to two referees whose detailed comments resulted in corrections and other improvements over an earlier version of this paper.

**2. Preliminaries; the partition function; strong disorder.** In this section we will first recall some basic facts about the partition function $Z_t$, and then give briefly some notions of Gaussian analysis which will be used later on. Let us recall that $W$ is a centered Gaussian field defined on $\mathbb{R}_+ \times \mathbb{R}$, which can also be seen as a Gaussian family $\{W(\varphi)\}$ indexed by tests functions $\varphi : \mathbb{R}_+ \times \mathbb{R} \to \mathbb{R}$, where $W(\varphi)$ stands for the Wiener integral of $\varphi$ with respect to $W$:

$$W(\varphi) = \int_\mathbb{R} \int_{\mathbb{R}_+} \varphi(s,x) W(ds,x)\, dx,$$

whose covariance structure is given by

$$(2.1) \quad \mathbf{E}[W(\varphi)W(\psi)] = \int_{\mathbb{R}_+} \left( \int_{\mathbb{R} \times \mathbb{R}} \varphi(s,x) Q(x-y) \psi(s,y) \, dx \, dy \right) ds,$$

for two arbitrary test functions $\varphi, \psi$.

Let us start here by defining more rigorously the quantity $H_t(b)$ given by (1.2), which can be done through a Fourier transform procedure: there exists (see, e.g., [5] for further details) a centered Gaussian independently scattered $\mathbb{C}$-valued measure $\nu$ on $\mathbb{R}_+ \times \mathbb{R}$ such that

$$(2.2) \quad W(t,x) = \int_{\mathbb{R}_+ \times \mathbb{R}} \mathbf{1}_{[0,t]}(s) e^{iux} \nu(ds, du).$$

For every test function $f : \mathbb{R}_+ \times \mathbb{R} \to \mathbb{C}$, set now

$$(2.3) \quad \nu(f) \equiv \int_{\mathbb{R}_+ \times \mathbb{R}} f(s,u) \nu(ds, du).$$

While the random variable $\nu(f)$ may be complex-valued, to ensure that it is real valued, it is sufficient to assume that $f$ is of the form $f(s,u) = f_1(s) e^{iu f_2(s)}$ for real valued functions $f_1$ and $f_2$. Then the law of $\nu$ is defined by the following covariance structure: for any such test functions $f, g : \mathbb{R}_+ \times \mathbb{R} \to \mathbb{C}$, we have

$$(2.4) \quad \mathbf{E}[\nu(f)\nu(g)] = \int_{\mathbb{R}_+ \times \mathbb{R}} f(s,u) \overline{g(s,u)} \hat{Q}(du) \, ds,$$

where the finite positive measure $\hat{Q}$ is the Fourier transform of $Q$ (see [17] for details).

From (2.2), we see that the Itô-stochastic differential of $W$ in time can be understood as $W(ds, x) := \int_{u \in \mathbb{R}} e^{iux} \nu(ds, du)$, or even, if the measure $\hat{Q}(du)$ has a density $f(u)$ with respect to the Lebesgue measure, which is typical, as

$$W(ds, x) := \int_{u \in \mathbb{R}} e^{iux} \sqrt{f(u)} M(ds, du),$$

where $M$ is a white-noise measure on $\mathbb{R}_+ \times \mathbb{R}$, that is, a centered independently scattered Gaussian measure with covariance given by $\mathbf{E}[M(A)M(B)] = m_{Leb}(A \cap B)$, where $m_{Leb}$ is Lebegue's measure on $\mathbb{R}_+ \times \mathbb{R}$.

We can go back now to the definition of $H_t(b)$: invoking the representation (2.2), we can write

$$(2.5) \quad -H_t(b) = \int_0^t W(ds, b_s) = \int_0^t \int_{\mathbb{R}} e^{iub_s} \nu(ds, du),$$

and it can be shown (see [5]) that the right-hand side of the above relation is well defined for any Hölder continuous path $b$, by a $L^2$-limit procedure.

Such a limiting procedure can be adapted to the specific case of constructing $H_t(b)$, using the natural time evolution structure; we will not comment on this further. However, the reader will surmise that the following remark, given for the sake of illustration, can be useful: when $\hat{Q}$ has a density $f$, we obtain

$$-H_t(b) = \iint_{[0,t]\times\mathbb{R}} e^{iub_s}\sqrt{f(u)}M(ds,du).$$

With the so-called *partition function* $Z_t^x$ defined earlier as $Z_t^x = E_b[e^{-\beta H_t(b)}]$, set

(2.6) $$p_t(\beta) := \frac{1}{t}\mathbf{E}[\log(Z_t^x)],$$

usually called the free energy of the system. By spatial homogeneity of $W$, $p_t(\beta)$ is independent of the initial condition $x \in \mathbb{R}$, and the same holds for the law of $b - x$ under $G_t^x$, thus without loss of generality, we set $x = 0$, hence, the notation $E_b, Z_t, \ldots$ standing for $E_b^0, Z_t^0$, etc. It was shown in [16] that $\lim_{t\to\infty} p_t(\beta) = \sup_{t\geq 0} p_t(\beta)$ exists and is positive, and that **P**-almost surely, $\frac{1}{t}\log Z_t$ converges to the same limit. The trivial bound

(2.7) $$p(\beta) := \lim_{t\to\infty} p_t(\beta) \leq \frac{\beta^2}{2}Q(0)$$

always holds, but the polymer is said to be in the *strong disorder* regime if $\lim_{t\to\infty}\frac{1}{t}\log Z_t < \frac{\beta^2}{2}Q(0)$, which is therefore equivalent to saying that inequality (2.7) above is strict. We will show in a separate publication that, for all $\beta \geq \beta_0$, while $p(\beta) \geq c\beta^{4/3}$ for all nontrivial random media $W$ and some constant depending on $W$'s law, we have the specific strong disorder upper bound $p(\beta) \leq c\beta^{2-2H/(2H+1)}$, where $H$ is a spatial Hölder exponent for $W$. Yet we do not know if these results can be made to hold for small $\beta$. One would prefer not having any condition on the temperature scale, and physicists expect strong disorder in our one-dimensional setting for all $\beta > 0$, which is only confirmed mathematically in some cases, such as in [3] and [6].

This is where the polymer's superdiffusivity (wandering exponent $\alpha > 1/2$) can be useful to our fully continuous situation. Since the concept of "strong disorder" was introduced in order to determine whether the random environment has any significant influence on polymer paths $b$, it is generally acceptable to say that a polymer with super-diffusive behavior exhibits "strong disorder." Even though this second definition does not match the common one given above $(p(\beta) = Q(0)\beta^2/2)$, it is useful to note that the results of the next section imply the following (see Corollary 7.3): if $W$ exhibits decorrelation that is not too slow, specifically if for large $x$, $Q(x) \leq cx^{-5/2-\vartheta}$, where $\vartheta > 0$, then the polymer is superdiffusive with exponent any $\alpha < \min\{\frac{1}{2} + \frac{\vartheta}{6-2\vartheta}; 3/5\}$, and this form of strong disorder holds for all $\beta > 0$. The specific order of decorrelation $x^{-5/2-\vartheta} \ll x^{-5/2}$ can be quantified by saying that $W$'s decorrelation is certainly faster than the well-known order $x^{-2+2H}$ for the increments of fractional Brownian motion, but the class of such $W$'s still qualifies as containing long-range correlations (polynomial with moderate power).

We also plan to investigate, in a separate publication, situations in which we can show the complementary story: we plan to prove that if *weak disorder* holds, that is, if $\lim_{t\to\infty} \frac{1}{t} \log Z_t = \frac{\beta^2}{2} Q(0)$, then the polymer is diffusive, that is, $\alpha = 1/2$.

**3. Discussion of hypothesis and results; strategy of proof.** Recall our goal: we will prove that, for the polymer measure $G_t = G_t^0$ in (1.3), Theorem 1.1 holds. This theorem gives an indication of the asymptotic speed of our polymer. Indeed, if we could write that $\sup_{s \leq t} |b_s| \sim t^\alpha$ under $G_t$ as $t \to \infty$, then Theorem 1.1 would state that the wandering exponent $\alpha$ is no smaller than $3/5$. As stated in the introduction, our basic technical assumption to prove the theorem is the following.

HYPOTHESIS 3.1. *We assume that $Q : \mathbb{R} \to \mathbb{R}$ defined by (1.1) is a symmetric positive function, decreasing on $\mathbb{R}_+$ and such that there exists a strictly positive constant $\theta$ such that*

$$Q(x) = O\left(\frac{1}{|x|^{3+\theta}}\right) \qquad as\ x \to \pm\infty.$$

The rate $3 + \theta$ can be quantified physically by saying that $W$ decorrelates in space faster than the well-known order $x^{-2+2H}$ for the increments of fractional Brownian motion with Hurst parameter $H \in (0,1)$, but the class of $W$'s defined by Hypothesis 3.1 still qualifies as containing long-range correlated noises (polynomial rate with moderate power), as opposed to exponential correlation decay, found, for instance, in finite memory ARCH/GARCH models, and even more so in opposition to the case of spatial white noise.

The specific correlation decay rate of $Q$ in the above hypothesis appears to be important in order to obtain the highest possible superdiffusion wandering exponent $\alpha$ using our technique (any $\alpha < 3/5$). The end of Section 7 shows that if one tries to use a smaller decay power than $3 + \theta$ above, the result is impeded: $\alpha$ cannot be chosen arbitrarily close to $3/5$. In Corollary 7.3 and its preceeding discussion, we prove that if $Q(x) = O(|x|^{-r})$ with $r \in (5/2, 3]$, then we can only guarantee being able to take $1/2 < \alpha < 3/(11 - 2r)$, so superdiffusivity is still proved, but $\alpha$ arbitrarily close to $3/5$ is disallowed.

Corollary 7.3 thus opens the interesting question of whether, in continuous space, the Brownian polymer in a Gaussian environment has a super-diffusive behavior with a wandering exponent determined by the environment's range/rate of spatial correlations. We do not believe that any physical conjecture in which $\alpha = 2/3$ specifically argues that this should hold in our continuous space setting. There are other examples in which scaling limits depend heavily on whether one is in discrete or continuous space: for instance, in the regime of small diffusion constant (resp. viscosity) $\kappa$, the almost-sure Lyapunov exponent for the partition function $Z_t$ (resp. Anderson model) is known to depend heavily on the spatial regularity of $W$ in continuous space (see [9]), but is known to be universally of order $1/\log(\kappa^{-1})$ in discrete space (see [4]). We will not discuss this point further herein.

REMARK 3.2. Hypothesis 3.1 immediately implies that $Q(0) < \infty$. Since $\max |Q| = Q(0)$ and $Q$ has an integrable tail, we get $Q \in L^1(\mathbb{R})$.

Without loss of generality, we will assume throughout that $Q$ is normalized so that $\int_{\mathbb{R}} Q(x)\,dx = 1$.

The integrability of $Q$ represents a kind of nondegeneracy condition, which says that the decorrelation of $W$ at distinct sites is not immediate.

STRATEGY OF THE PROOF FOR THEOREM 1.1. For $t, \epsilon > 0$, set

$$A_{t,\epsilon} = \{\text{there exists } s_0 \in [t/2, t] \text{ such that } |b_{s_0}| \geq t^{3/5 - \epsilon/2}\}.$$

Then we can write

$$\frac{\langle \sup_{s \leq t} |b_s| \rangle_t}{t^{3/5-\epsilon}} \geq \frac{t^{\epsilon/2}}{t^{3/5-\epsilon/2}} \left\langle \sup_{s \leq t} |b_s| \mathbf{1}_{A_{t,\epsilon}} \right\rangle_t \geq t^{\epsilon/2} G_t(A_{t,\epsilon}),$$

since $\sup_{s \leq t} |b_s| \geq t^{3/5-\epsilon/2}$ on $A_{t,\epsilon}$. Thus,

(3.1) $$\frac{\langle \sup_{s \leq t} |b_s| \rangle_t}{t^{3/5-\epsilon}} \geq t^{\epsilon/2}(1 - G_t(A_{t,\epsilon}^c)),$$

where $A_{t,\epsilon}^c = \{b; \sup_{s \in [t/2,t]} |b_s| \leq t^{3/5-\epsilon/2}\}$ is the complement of $A_{t,\epsilon}$. We will start now a discretization procedure in space: for an arbitrary integer $k$, and $\alpha > 0$, set

$$I_k^\alpha = t^\alpha [2k-1, 2k+1) \quad \text{and} \quad L_k^\alpha = \{b; b_s \in I_k^\alpha \text{ for all } s \in [t/2, t]\}.$$

Then $\tilde{A}_{t,\epsilon} = L_0^{3/5-\epsilon/2}$, and equation (3.1) can be rewritten as

$$\frac{\langle \sup_{s \leq t} |b_s| \rangle_t}{t^{3/5-\epsilon}} \geq t^{\epsilon/2}(1 - G_t(L_0^{3/5-\epsilon/2})).$$

Set now

$$Z_t^\alpha(k) := E_b[\mathbf{1}_{L_k^\alpha} \exp(-\beta H_t(b))].$$

We have
$$\frac{\langle \sup_{s\leq t}|b_s|\rangle_t}{t^{3/5-\epsilon}} \geq t^{\epsilon/2}\left(1 - \frac{Z_t^{3/5-\epsilon/2}(0)}{E_b[\exp(-\beta H_t(b))]}\right),$$
by definition of $G_t$. On the other hand, since the events $L_k^\alpha$ are disjoint sets, we have
$$E_b[\exp(-\beta H_t(b))] \geq \sum_{k\in\mathbb{Z}} Z_t^{3/5-\epsilon/2}(k).$$
Therefore, we have established that
$$(3.2)\qquad \frac{\langle \sup_{s\leq t}|b_s|\rangle_t}{t^{3/5-\epsilon}} \geq t^{\epsilon/2}\left(1 - \frac{Z_t^{3/5-\epsilon/2}(0)}{Z_t^{3/5-\epsilon/2}(0) + Z_t^{3/5-\epsilon/2}(k)}\right),$$
for any integer $k \neq 0$. Suppose now that $W \in \mathcal{A}_t$, where $\mathcal{A}_t$ is defined as
$$\mathcal{A}_t := \{W; \text{There exists } k^* \neq 0 \text{ such that } Z_t^\alpha(k^*) > Z_t^\alpha(0)\}.$$
Then, choosing $k = k^*$ in (3.2), it is easily seen that
$$\frac{\langle \sup_{s\leq t}|b_s|\rangle_t}{t^{3/5-\epsilon}} \geq t^{\epsilon/2}\left(1 - \frac{1}{2}\right) \geq 1,$$
whenever $t$ is large enough. The proof is now easily finished if we can prove the following lemma:

LEMMA 3.3. *Given a positive real number $\alpha \in (1/2; 3/5)$ and an environment $W$ satisfying Hypothesis 3.1, then*
$$(3.3)\qquad \liminf_{t\to\infty} \mathbf{P}(\mathcal{A}_t) = 1.$$

The remainder of this article will now be devoted to the proof Lemma 3.3. $\square$

**4. Initial covariance computations.** In order to prove Lemma 3.3, we shall begin with a series of preliminary results, the first of which is a covariance computation, including precise asymptotic estimations in large time, for space-time averages of the random environment $W$.

For a given $k \in \mathbb{Z}$ and $\alpha > 0$, recall that $I_k := I_k^\alpha = t^\alpha[2k-1, 2k+1)$, and set
$$(4.1)\qquad \tilde{\eta}_k = \tilde{\eta}_k^\alpha := \frac{1}{t^{(\alpha+1)/2}} \int_{t/2}^t \int_{I_k} W(ds, x)\, dx.$$
Then $\{\tilde{\eta}_k; k \in \mathbb{Z}\}$ is a centered Gaussian vector, whose covariance matrix will be called $C(t) = (C_{\ell,k}(t))_{\ell,k\in\mathbb{Z}}$, where
$$(4.2)\qquad C_{\ell,k}(t) = \mathbf{E}[\tilde{\eta}_\ell \tilde{\eta}_k] = \mathbf{Cov}(\tilde{\eta}_\ell; \tilde{\eta}_k) = \frac{1}{2t^\alpha} \int_{I_k} \int_{I_\ell} Q(x-y)\, dx\, dy,$$

where the last equality above follows directly from the definition of $W$'s covariance in (2.1).

Here and below, we omit the superscripts $\alpha$ on quantities like $\tilde{\eta}_k^\alpha$, $I_k^\alpha$, $L_k^\alpha$, etc. We now proceed to estimate the matrix $C(t)$, and show, in particular, that $\lim_{t \to \infty} C(t) = \text{Id}$. This can be interpreted as saying that the amount of decorrelation of the potential at distant locations implied by Hypothesis 3.1 is enough to guarantee independence of the $\tilde{\eta}_k$ asymptotically.

PROPOSITION 4.1. *Let $\theta$ be the strictly positive constant defined in Hypothesis 3.1, and consider $k \in \mathbb{Z}$, $\alpha > 0$ and $\tau < \theta \wedge 1$. Set also*

$$\lambda := \frac{1}{C_{0,0}(t)} = \frac{1}{C_{k,k}(t)},$$

*where $C(t)$ has been defined at (4.2). Then, the elements of $C(t)$ satisfy the following:*

(i) $\lambda = 1 + O(\frac{1}{t^\alpha})$.
(ii) $\lambda \sum_{\ell \neq k} |\ell - k|^\tau |C_{\ell,k}(t)| = O(\frac{1}{t^\alpha})$.

PROOF.

*Step* 0: *initial calculation.* We will only consider the case $k = 0$, the other ones being easily deduced by homogeneity of $W$. Let us first evaluate $C_{\ell,0}(t)$ for $\ell \geq 0$ (here again, the case $\ell < 0$ is similar, since $Q$ is a symmetric function). Then, a direct application of (4.2) gives

$$C_{\ell,0}(t) = \frac{1}{2t^\alpha} \int_{t^\alpha(2\ell-1)}^{t^\alpha(2\ell+1)} \int_{-t^\alpha}^{t^\alpha} Q(x-y)\,dx\,dy.$$

Set now

$$(I) := \frac{1}{2t^\alpha} \left[ \int_{t^\alpha(2\ell-1)}^{t^\alpha(2\ell+1)} \int_{-\infty}^{-t^\alpha} Q(x-y)\,dx\,dy + \int_{t^\alpha(2\ell-1)}^{t^\alpha(2\ell+1)} \int_{t^\alpha}^{\infty} Q(x-y)\,dx\,dy \right].$$

Since $\int_{\mathbb{R}} Q(x-y)\,dx = 1$ for any $y \in \mathbb{R}$, it is easily checked that

(4.3) $$C_{\ell,0}(t) = 1 - (I).$$

Then, a series of changes of variable yields

$$(I) = \frac{1}{2t^\alpha} \left[ \int_{t^\alpha(2\ell-1)}^{t^\alpha(2\ell+1)} \int_{-\infty}^{-t^\alpha-y} Q(u)\,du\,dy + \int_{t^\alpha(2\ell-1)}^{t^\alpha(2\ell+1)} \int_{t^\alpha-y}^{\infty} Q(u)\,du\,dy \right]$$

$$= \frac{1}{2t^\alpha} \left[ \int_{-t^\alpha(2\ell+2)}^{-t^\alpha(2\ell)} \int_{-\infty}^{\hat{z}} Q(u)\,du\,d\hat{z} + \int_{-t^\alpha(2\ell)}^{-t^\alpha(2\ell-2)} \int_{z}^{\infty} Q(u)\,du\,dz \right],$$

where we have set $\hat{z} = -t^\alpha - y$ and $z = t^\alpha - y$. Thus, denoting by $\bar{F}(z)$ the quantity $\int_z^\infty Q(u)\,du$, we get

$$
\begin{aligned}
(I) &= \frac{1}{2t^\alpha}\left[\int_{-t^\alpha(2\ell+2)}^{-t^\alpha(2\ell)}(1-\bar{F}(\hat{z}))\,d\hat{z} + \int_{-t^\alpha(2\ell)}^{-t^\alpha(2\ell-2)} \bar{F}(z)\,dz\right] \\
&= 1 - \frac{1}{2t^\alpha}\int_{-t^{-\alpha}(2\ell+2)}^{-t^\alpha(2\ell)}\bar{F}(z)\,dz + \frac{1}{2t^\alpha}\int_{-t^\alpha(2\ell)}^{-t^\alpha(2\ell-2)}\bar{F}(z)\,dz.
\end{aligned}
\tag{4.4}
$$

Putting together (4.3) and (4.4), one obtains, for any $\ell \geq 0$,

$$
C_{\ell,0}(t) = \frac{1}{2t^\alpha}\left[\int_{-t^\alpha(2\ell+2)}^{-t^\alpha(2\ell)}\bar{F}(z)\,dz - \int_{-t^\alpha(2\ell)}^{-t^\alpha(2\ell-2)}\bar{F}(z)\,dz\right].
\tag{4.5}
$$

*Step* 1: *proving item* (i). We are now ready to prove item (i). By symmetry of $Q$, we have $1 - \bar{F}(-z) = \bar{F}(z)$. Thus, for $\ell = 0$, equation (4.5) becomes

$$
\begin{aligned}
C_{0,0}(t) &= \frac{1}{2t^\alpha}\left[\int_{-2t^\alpha}^{0}(1-\bar{F}(-z))\,dz - \int_0^{2t^\alpha}\bar{F}(z)\,dz\right] \\
&= 1 - \frac{1}{t^\alpha}\int_0^{2t^\alpha}\bar{F}(z)\,dz.
\end{aligned}
\tag{4.6}
$$

Now, using the fact that

$$
\bar{F}(z) \leq c(1 \wedge |z|^{-(2+\theta)}),
\tag{4.7}
$$

which follows directly from Hypothesis 3.1, it is easily seen that $C_{0,0}(t) = 1 + O(t^{-\alpha})$, which ends the proof of item (i).

*Step* 2: *proving item* (ii). In order to show item (ii), we deal with $\ell = 1$ separately from the other cases. Beginning with $\ell \geq 2$, we first get the obvious derivative $\bar{F}'(z) = -Q(z)$, and we will use the fact that $Q$ is decreasing on $\mathbb{R}_+$ to bound this latter function on an interval in $\mathbb{R}_+$ by its value at the left endpoint. Invoking the fact that $\bar{F}(-v) = 1 - \bar{F}(x)$, we may thus write from equation (4.5)

$$
\begin{aligned}
|C_{\ell,0}(t)| &= \frac{1}{2t^\alpha}\left|\int_{-t^\alpha(2\ell)}^{-t^\alpha(2\ell-2)}[\bar{F}(z-2t^\alpha) - \bar{F}(z)]\,dz\right| \\
&= \frac{1}{2t^\alpha}\left|\int_{-t^\alpha(2\ell)}^{-t^\alpha(2\ell-2)}[\bar{F}(-z+2t^\alpha) - \bar{F}(-z)]\,dz\right| \\
&= \frac{1}{2t^\alpha}\left|\int_{t^\alpha(2\ell-2)}^{2t^\alpha\ell}[\bar{F}(z+2t^\alpha) - \bar{F}(z)]\,dz\right| \\
&= \frac{1}{2t^\alpha}\left|\int_{t^\alpha(2\ell-2)}^{2t^\alpha\ell}\left(-\int_z^{z+2t^\alpha}Q(x)\,dx\right)dz\right| \\
&\leq 2t^\alpha Q(t^\alpha(2\ell-2)) \\
&\leq c t^{-\alpha(2+\theta)}(2\ell-2)^{-3-\theta},
\end{aligned}
$$

where the last step holds by Hypothesis 3.1 for some constant $c > 0$. We immediately obtain

$$\sum_{\ell=2}^{\infty} |C_{\ell,0}(t)| \ell^\tau \leq ct^{-\alpha(2+\theta)} \sum_{\ell=2}^{\infty} (2\ell - 2)^{-3-\theta} \ell^\tau$$

$$\leq cK_{\tau,\theta} t^{-\alpha(2+\theta)}$$

for some constant $K_{\tau,\theta}$ as soon as $\tau < 2 + \theta$, which is clearly satisfied by the assumption on $\tau$, and leads to an upper bound in the series in item (ii) which is amply sufficient to prove the proposition, except for the term $\ell = 1$, with which we deal now.

To finish the proof of the proposition, it is indeed sufficient to prove that $t^\alpha C_{1,0}$ is bounded. We first evaluate this quantity from (4.5):

$$t^\alpha C_{1,0} = \int_{-4t^\alpha}^{-2t^\alpha} \bar{F}(z)\,dz - \int_{-2t^\alpha}^{0} \bar{F}(z)\,dz$$

$$= \int_{-2t^\alpha}^{0} (\bar{F}(z - 2t^\alpha) - \bar{F}(z))\,dz$$

$$= \int_{-2t^\alpha}^{0} \left( \int_{z-2t^\alpha}^{z} Q(x)\,dx \right) dz$$

$$= \int_{0}^{2t^\alpha} \left( \int_{-z-2t^\alpha}^{-z} Q(x)\,dx \right) dz$$

$$= \int_{0}^{2t^\alpha} \left( \int_{z}^{z+2t^\alpha} Q(x)\,dx \right) dz.$$

Next we separate the first unit of the $z$-integral from its remainder: $t^\alpha C_{1,0} = A + B$, where we define $A := \int_0^1 (\int_z^{z+2t^\alpha} Q(x)\,dx)\,dz$ and $B := \int_1^{1 \wedge 2t^\alpha}(\int_z^{z+2t^\alpha} Q(x)\,dx)\,dz$. Since $\int_\mathbb{R} Q = 1$, we immediately have $A \leq 1$ which is the only term to deal with when $t \leq 2^{-1/\alpha}$. When $t > 2^{-1/\alpha}$, for the term $B$, we use Hypothesis 3.1: for some constant $c$,

$$B \leq c \int_1^{2t^\alpha} \left( \int_z^{z+2t^\alpha} x^{-3-\theta}\,dx \right) dz$$

$$= \frac{c}{(\theta+1)(\theta+2)}(1 - 2^{-\theta} + 4^{-\theta-1})(t^\alpha)^{-(\theta+1)} \leq \frac{c}{(\theta+1)(\theta+2)}.$$

This finishes the proof of the proposition. □

**5. Interaction between $b$ and $W$.** The next step in developing the tools to prove Lemma 3.3 is to get some quantitative information about the way $b$ interacts with the random environment $W$ when the Brownian motion is localized by the event $L_k$. As we did with the notation $I_k := I_k^\alpha$, we are omitting superscripts $\alpha$, writing only $L_k$ instead of $L_k^\alpha$ from now on.

We begin by introducing two quantities. First, in order to simplify some $t$-dependent normalizers, we renormalize $\tilde{\eta}$ as

$$(5.1) \qquad \eta_\ell := \frac{t^{(1-\alpha)/2}}{2}\tilde{\eta}_\ell = \frac{1}{2t^\alpha}\int_{t/2}^t \int_{I_\ell} W(ds,x)\,dx;$$

we will not need to revert to using $\tilde{\eta}$ in this article. We also need a vector $v = v(b_s; t/2 \le s \le t)$ of $\mathbb{R}^\mathbb{Z}$, defined for each $\ell \in \mathbb{Z}$ by

$$(5.2) \qquad v_\ell := 4t^{\alpha-1}\mathbf{E}\left[\eta_\ell \int_{t/2}^t W(ds, b_s)\right].$$

We will prove, in a sense, that $v_\ell$ looks like $\mathbf{1}_{\{k\}}(\ell)$ on $L_k$. To this end, for a fixed $k \in \mathbb{Z}$, and $\tau < \theta$ (remember that $\theta$ is defined in Hypothesis 3.1), let us consider the norm $\|\cdot\|_{\tau,k}$ defined on $\mathbb{R}^\mathbb{Z}$ by

$$(5.3) \qquad \|x\|_{\tau,k} = |x_k| + \sum_{i \ne k} |x_i||i-k|^\tau.$$

REMARK 5.1. It will be essential in the sequel to control the decay of $v_\ell$, and also of a quantity $\delta_\ell$ [defined later in Proposition 5.3 as the $\ell$th component of the solution $x$ to the linear system $C(t)x = v$] when $|\ell| \to \infty$. It will be used, for instance, in relations (7.6) and (7.10). This is why we have introduced the norm $\|\cdot\|_{\tau,k}$ here.

5.1. *Asymptotics and boundedness of $v$.* We are now ready to state a first result about the interaction between $b$ and $W$: the behavior of the vector $v$ in large time.

PROPOSITION 5.2. *Suppose $b \in L_k$. Then the vector $v$ given by (5.2) satisfies the following properties:*

(i) *Let $\|\cdot\|_{\tau,k}$ be the norm defined at (5.3). Then*

$$\|v\|_{\tau,k} - v_k = O\left(\frac{1}{t^\alpha}\right).$$

(ii) *For $t$ large enough, there exist two strictly positive real numbers $\underline{c}$ and $\overline{c}$ such that*

$$\underline{c} \le v_k \le \overline{c}.$$

PROOF. Let us start with item (i). To perform calculations rigorously, it is best to use the environment representation (2.2). Recall also that $\eta_k$ is given by (5.1). Then

$$v_\ell = \frac{2}{t}\mathbf{E}\left[\int_{t/2}^{t}\int_{\mathbb{R}}\exp(iub_s)\nu(ds,du)\int_{I_\ell}\int_{t/2}^{t}\int_{\mathbb{R}}\exp(iux)\nu(ds,du)\,dx\right]$$

$$= \frac{2}{t}\int_{I_\ell}\mathbf{E}\left[\int_{t/2}^{t}\int_{\mathbb{R}}\exp(iub_s)\nu(ds,du)\int_{t/2}^{t}\int_{\mathbb{R}}\exp(iux)\nu(ds,du)\right]dx.$$

Thanks to (2.4), and according to the fact that $\hat{Q}$ is the Fourier transform of $Q$, we thus have

$$v_\ell = \frac{2}{t}\int_{I_\ell}\left[\int_{t/2}^{t}\int_{\mathbb{R}}\exp(iu(b_s-x))\hat{Q}(du)\,ds\right]dx$$

(5.4)
$$= \frac{2}{t}\int_{t/2}^{t}\int_{I_\ell}Q(b_s-x)\,dx\,ds$$

(5.5)
$$\leq \sup_{s\in[t/2,t]}\int_{I_\ell}Q(b_s-x)\,dx.$$

However, if $\ell \neq k$, on the event $L_k$, it is easily checked that, for $s \in [t/2,t]$, and for all $x \in I_\ell$, we have

$$(2|\ell-k|-2)t^\alpha \leq |b_s-x|.$$

According to the fact that $Q$ is a positive decreasing function on $\mathbb{R}_+$, and $Q(x) = Q(|x|)$, for each $s \in [t/2,t]$ we can conclude that

$$\int_{I_\ell}Q(b_s-x)\,dx = \int_{I_\ell}Q(|b_s-x|)\,dx$$

$$\leq \int_{t^\alpha(2\ell+1)}^{t^\alpha(2\ell+1)}Q((2|\ell-k|-2)t^\alpha)\,dx$$

$$\leq 2t^\alpha Q(t^\alpha(2|\ell-k|-2)).$$

Consequently, putting together equations (5.5) and (5.1), we get

$$\|v\|_{\tau,k} = v_k + \sum_{\ell\neq k}|\ell-k|^\tau v_\ell$$

$$\leq v_k + 2t^\alpha\sum_{\ell\neq k}|\ell-k|^\tau Q(t^\alpha(2|\ell-k|-2))$$

(5.6)
$$\leq v_k + \frac{\kappa}{t^{\alpha(2+\theta)}}\sum_{\ell\neq k}|\ell-k|^{-(3+\theta-\tau)}$$

$$\leq v_k + \frac{\kappa}{t^{\alpha(2+\theta)}},$$

where $\kappa$ is a positive constant that can change from one occurrence to the next, and where we have used again Hypothesis 3.1. It is now readily checked that $\|v\|_{\tau,k} \leq v_k + O(t^{-\alpha})$, which ends the proof of item (i).

Let us prove now item (ii): go back to equation (5.4) and set $\ell = k$. Then we get

$$\inf_{s \in [t/2,t]} \int_{I_k} Q(b_s - x)\, dx \leq v_k \leq \sup_{s \in [t/2,t]} \int_{I_k} Q(b_s - x)\, dx$$

$$\leq \int_{\mathbb{R}} Q(u)\, du = 1.$$

To find a lower bound on the left-hand side, we now make use of the nondegeneracy assumption, as noted in Remark 3.2: since $Q$ is an even function, we get $\int_0^\infty Q(x)\, dx = 1/2$. But if $b \in L_k$, then for any $s \in [t/2, t]$, we have that the interval $b_s - I_k$ contains either $[0, t^\alpha]$ or $[-t^\alpha, 0]$, so that, again by the evenness of $Q$,

$$\int_{I_k} Q(b_s - x)\, dx \geq \int_0^{t^\alpha} Q(x)\, dx.$$

The latter quantity, which tends to $1/2$ when $t \to \infty$, can be made to exceed $1/4$ for $t$ large enough. This finishes the proof of item (ii) with $\underline{c} = 1/4$ and $\overline{c} = 1$, and the proposition. $\square$

5.2. *Inversion of $C(t)$.* In this section we will be concerned with the operator $C^{-1}(t)$, where $C(t)$ has been defined by relation (4.2), and more specifically, we will get some information about the solution $\delta$ to the system $C(t)x = v$. The importance of $\delta$ stems from the fact that the variables $\eta_k$ will be independent of $-H_t(b) - \sum_{j \in \mathbb{Z}} \delta_j \eta_j$, which will be useful for further computations (see Proposition 7.2). However, we have already seen that $C(t)$ behaves asymptotically like the identity matrix, and thus, the vector $\delta$ should be of the same kind as $v$, in particular, when $b \in L_k$. This is indeed the case, and will be proved in the following proposition.

PROPOSITION 5.3. *Under Hypothesis 3.1, suppose in addition that $b \in L_k$. Set $l_{\tau,k} = \{x \in \mathbb{R}^{\mathbb{Z}}; \|x\|_{\tau,k} < \infty\}$. Then:*

(i) *The operator $C(t)$ is invertible in $l_{\tau,k}$. We set then $\delta := C^{-1}(t)v$.*
(ii) *There exist some strictly positive real numbers $\underline{d}$ and $\overline{d}$ such that*

$$\underline{d} \leq \delta_k \leq \overline{d}.$$

(iii) *The following relation holds:*

$$\|\delta\|_{\tau,k} - \delta_k = O\left(\frac{1}{t^\alpha}\right).$$

(iv) *On the probability space* $(\Omega, \mathcal{G}, \mathbf{P})$, *the family* $\{\eta_l; l \in \mathbb{Z}\}$ *is independent of* $-H_t(b) - \sum_{j \in \mathbb{Z}} \delta_j \eta_j$.

REMARK 5.4. Notice that Proposition 5.3 contains a considerable amount of the information which will be used for the proof of Lemma 3.3. Indeed, inequality (7.14) will be obtained thanks to item (iv), item (iii) will be invoked for inequality (7.10), and item (ii) will be essential in order to define the random variables $\breve{\eta}_0$ and $\breve{\eta}_k$ in (7.9).

PROOF OF PROPOSITION 5.3.

*Step* 1: *proving item* (i). We choose the standard operator norm on $l_{\tau,k}$: a matrix $A$ is defined to be in the linear operator space $\mathcal{L}_{\tau,k}$ if the norm

$$\|A\|_{\tau,k} := \sup_{x \in l_{\tau,k}: \|x\|_{\tau,k}=1} \|Ax\|_{\tau,k}$$

is finite. Then, on one hand, the following relations are satisfied since we are dealing with the operator norm on $l_{\tau,k}$: for $D_1, D_2 \in \mathcal{L}_{\tau,k}$ and $x \in l_{\tau,k}$:

(5.7) $\quad \|D_1 x\|_{\tau,k} \leq \|D_1\|_{\tau,k} \|x\|_{\tau,k}$ and $\|D_1 + D_2\|_{\tau,k} \leq \|D_1\|_{\tau,k} + \|D_2\|_{\tau,k}$.

On the other hand, let us now prove that, setting $A(t) := \text{Id} - \lambda C(t)$, Proposition 4.1 yields that $\|A(t)\|_{\tau,k} = O(t^{-\alpha})$, and thus,

(5.8) $$\|A(t)\|_{\tau,k} < 1,$$

if $t$ is large enough. First recall that by definition of $C(t)$ and $\lambda$, denoting by $\dot{C}(t)$ the matrix $C(t)$ deprived of its diagonal, we have

$$A(t) = -\lambda \dot{C}(t).$$

By Proposition 4.1 item (i), $\lambda$ tends to 1 as $t \to \infty$. Therefore, it is sufficient to show that $\|\dot{C}(t)\|_{\tau,k} = O(t^{-\alpha})$. Thus, let $x \in l_{\tau,k}$ such that $\|x\|_{\tau,k} = 1$. In other words,

$$|x_k| + \sum_{i \neq k} |x_i||i - k|^\tau = 1.$$

Now we calculate the two terms that form $\|\dot{C}(t)x\|_{\tau,k}$. The first is

(5.9)
$$\begin{aligned}
|(\dot{C}(t)x)_k| &= \left|\sum_{j \neq k} C_{kj}(t) x_j\right| \leq \sum_{j \neq k} |C_{kj}(t) x_j| \\
&\leq \left(\sum_{j \neq k} |x_j||k-j|^\tau\right)\left(\sum_{j \neq k} |C_{kj}(t)||k-j|^\tau\right) \\
&\leq 1 \cdot O(t^{-\alpha}),
\end{aligned}$$

where we used the assumption $\|x\|_{\tau,k} = 1$ and the result of Proposition 4.1 item (ii). The second term in $\|\dot{C}(t)x\|_{\tau,k}$ equals

$$\sum_{i \neq k}\left|\sum_{j \neq i} C_{ij}(t)x_j\right||i-k|^\tau \leq \sum_{j \in \mathbb{Z}} |x_j| \sum_{i \neq j; i \neq k} |C_{ij}(t)||i-k|^\tau =: K_2;$$

we split this sum up according to $j = k$ or $j \neq k$:

$$K_2 \leq |x_k| \sum_{i \neq k} |C_{ik}(t)||i-k|^\tau + \sum_{j \neq k} |x_j| \sum_{i \neq j; i \neq k} |C_{ij}(t)||i - j + j - k|^\tau$$

$$\leq |x_k| \sum_{i \neq k} |C_{ik}(t)||i-k|^\tau + \sum_{j \neq k} |x_j| \sum_{i \neq j; i \neq k} |C_{ij}(t)||i-j|^\tau$$

$$+ \sum_{j \neq k} |x_j||j-k|^\tau \sum_{i \neq j; i \neq k} |C_{ij}(t)|,$$

where in the last line we used the fact that $|a+b|^\tau \leq |a|^\tau + |b|^\tau$ whenever $\tau \in (0, 1)$.

Now using the fact that $\sum_{i \neq j; i \neq k} |C_{ij}(t)|$ is bounded above by $\sum_{i \neq j} |C_{ij}(t)||i-j|^\tau$, and the latter is $O(t^{-\alpha})$ by Proposition 4.1 item (ii), we can assert $K_2 \leq O(t^{-\alpha})$, which, combined with (5.9), implies our goal $\|\dot{C}(t)\|_{\tau,k} = O(t^{-\alpha})$, and thus (5.8). This contraction relation (5.8) finishes the proof of (i) because it allows us to define $C^{-1}(t)$ in $\mathcal{L}_{\tau,k}$ by a Von Neumann type series of the form

$$(5.10) \qquad C^{-1}(t) = \lambda \sum_{j \geq 0} A^j.$$

*Step* 2: *proving item* (ii). For $t$ large enough, set $\delta = C^{-1}(t)v$, which makes sense since $v \in l_{\tau,k}$. Then, thanks to the fact that $C^{-1}(t)$ can be defined by relation (5.10), we have

$$\delta_k = \lambda\left(v_k + \sum_{j \geq 1}(A^j v)_k\right) \geq \lambda\left(v_k - \sum_{j \geq 1} \|A^j v\|_{\tau,k}\right)$$

$$\geq \lambda\left(v_k - \sum_{j \geq 1} \|A\|_{\tau,k}^j \|v\|_{\tau,k}\right),$$

where we have used the relations $x_k \geq -\|x\|_{\tau,k}$ and (5.7). Hence, since $\|A(t)\|_{\tau,k} = O(t^{-\alpha})$, we obtain

$$(5.11) \quad \delta_k \geq \lambda\left(v_k - \frac{\|A\|_{\tau,k}}{1 - \|A\|_{\tau,k}} \|v\|_{\tau,k}\right) \geq \lambda\left(v_k + O\left(\frac{1}{t^\alpha}\right)\right) \geq \underline{d} + O\left(\frac{1}{t^\alpha}\right),$$

according to the properties of $v$ shown at Proposition 5.2. The upper bound on $\delta_k$ can now be shown by the same type of argument, which ends the proof of our claim.

*Step* 3: *proving item* (iii). Let us evaluate now the quantity $\|\delta\|_{\tau,k} - \delta_k$: thanks to relations (5.7) and (5.11), we get

$$\|\delta\|_{\tau,k} - \delta_k \leq \|C(t)^{-1}\|_{\tau,k} \|v\|_{\tau,k} - \delta_k$$
$$\leq \left( \|C(t)^{-1}\|_{\tau,k} \|v\|_{\tau,k} - \lambda v_k + \frac{\lambda \|A\|_{\tau,k}}{1 - \|A\|_{\tau,k}} \|v\|_{\tau,k} \right).$$

Thus, using again the fact that $C^{-1}(t)$ is defined by equation (5.10) and relation (5.7), we obtain

$$\|\delta\|_{\tau,k} - \delta_k \leq \lambda \left( \frac{1 + \|A\|_{\tau,k}}{1 - \|A\|_{\tau,k}} \|v\|_{\tau,k} - v_k \right)$$
$$= \lambda (\|v\|_{\tau,k} - v_k) + O\left(\frac{1}{t^\alpha}\right) = O\left(\frac{1}{t^\alpha}\right),$$

where in the last two steps we have invoked, respectively, item (i) and Proposition 5.2. This concludes our proof of (iii).

*Step* 4: *proving item* (iv). Recall that, by definition, $C(t) = t^{-(1-\alpha)} \mathbf{Cov}(\eta)$. Hence,

$$\delta_j = (C^{-1}(t)v)_j = \tfrac{1}{4} t^{1-\alpha} \sum_{k \in \mathbb{Z}} [\mathbf{Cov}(\eta)]_{jk}^{-1} v_k$$
$$= \sum_{k \in \mathbb{Z}} [\mathbf{Cov}(\eta)]_{jk}^{-1} \mathbf{E}\left[ \int_{t/2}^{t} W(ds, b_s) \eta_k \right]$$
$$= \sum_{k \in \mathbb{Z}} [\mathbf{Cov}(\eta)]_{jk}^{-1} \mathbf{E}[(-H_t(b)) \eta_k];$$

we have the following standard calculation for any $\ell \in \mathbb{Z}$:

$$\mathbf{E}\left[ \left( -H_t(b) - \sum_{j \in \mathbb{Z}} \delta_j \eta_j \right) \eta_\ell \right]$$
$$= -\mathbf{E}[H_t(b) \eta_\ell] + \mathbf{E} \sum_{j \in \mathbb{Z}} \sum_{k \in \mathbb{Z}} [\mathbf{Cov}(\eta)]_{jk}^{-1} \mathbf{E}[H_t(b) \eta_k] \eta_j \eta_\ell$$
$$= -\mathbf{E}[H_t(b) \eta_\ell] + \mathbf{E} \sum_{j \in \mathbb{Z}} \sum_{k \in \mathbb{Z}} [\mathbf{Cov}(\eta)]_{jk}^{-1} [\mathbf{Cov}(\eta)]_{j\ell} \mathbf{E}[H_t(b) \eta_k]$$
$$= -\mathbf{E}[H_t(b) \eta_\ell] + \sum_{k \in \mathbb{Z}} \delta_{k\ell} \mathbf{E}[H_t(b) \eta_k] = 0.$$

Now since for fixed $b$, $H_t(b)$ and the sequence $\eta$ are both linear functionals of a same Gaussian field, they form a jointly Gaussian vector, and are thus independent. $\square$

**6. Application of Girsanov's theorem.** In our context the cost of having $b$ living in the interval $I_k = [t^\alpha(2k-1), t^\alpha(2k+1)]$ instead of $I_0 = [-t^\alpha, t^\alpha]$ can be calculated explicitly thanks to Girsanov's theorem: given an integer $k$, a real number $t$ and a realization of the environment $W$, we define a new environment by setting $W^{k,t}(ds, x) := W(ds, x + h(s))$, where

$$h(s) := \min(2s/t, 1) 2kt^\alpha,$$

or more rigorously,

$$(6.1) \qquad W^{k,t}(s, x) := \int_0^s W(du, x + h(u)).$$

A simple and useful result that we can now prove is the following.

LEMMA 6.1. *The random fields defined by $W = \{W(s, x) : (s, x) \in \mathbb{R}_+ \times \mathbb{R}\}$ and $W^{k,t} = \{\int_0^s W(du, u + h(u)) : (s, x) \in \mathbb{R}_+ \times \mathbb{R}\}$ have the same distribution.*

PROOF. The easiest way to establish this result is to revert to the representation of $W$ using the Gaussian measure $\nu$, that is, (2.2), and also its consequence (2.5), so that

$$W^{k,t}(s, x) := \int_0^s \int_\mathbb{R} e^{\imath\lambda(x + h(u))} \nu(ds, du).$$

Since the law of this centered Gaussian field is determined by its covariance structure only, it is now immediate to check, using the formulas (2.3) and (2.4), that it has the same law as $W$, since we have

$$W(s, x) := \int_0^s \int_\mathbb{R} e^{\imath\lambda x} \nu(ds, du).$$

The calculations are left to the reader. $\square$

ALTERNATE PROOF. It is also possible to invoke a direct proof of this fact, using $L^2$ approximations of $W^{k,t}(s, x)$ by Riemann sums. For fixed $s, x$, $W^{k,t}(s, x)$ can be written as a limit in $L^2(\Omega)$, as $n \to \infty$, of the sum $\sum_{i=1}^n J_i^{k,t}$ of the increments $J_i^{k,t} := W([si/n, s(i+1)/n], x + h(si/n))$, whose individual laws are identical to those of the $J_i$'s defined without adding the shift $h(si/n)$, because $W$ is spatially homogeneous. Since the $J_i^{k,t}$'s are independent as $i$ changes (as are the $J_i$'s), $W^{k,t}(s, x)$ and $W(s, x)$ have the same distribution for fixed $s, x$; we omit the end of this—more intuitive but less rigorous—proof. $\square$

We also need to introduce a modified partition function $\tilde{Z}$ defined by

$$(6.2) \qquad \tilde{Z}_t^\alpha(k) = E_b\left[\mathbf{1}_{L_k}(b) \exp\left(\beta\left(\int_0^t W(ds, b_s) - \sum_{j \in \mathbb{Z}} \delta_j \eta_j\right)\right)\right].$$

In the sequel we will have to stress the dependence of these partition functions on the environment under consideration. We will thus set $\tilde{Z}_t^\alpha(k) = \tilde{Z}_t^\alpha(k, W)$. With these notations in mind, we can prove the following proposition, which shows that the cost of having $b$ live in $L_k$ rather than $L_0$ is exponential of order $t^{2\alpha-1}$.

PROPOSITION 6.2. *Given two positive real numbers $\alpha$ and $t$, and an integer $k$ fixed, we have*

(6.3) $$\tilde{Z}_t^\alpha(k, W) \geq \exp[-4(k+k^2)t^{2\alpha-1}]\tilde{Z}_t^\alpha(0, W^{k,t}).$$

PROOF.
*Step* 1: *using Girsanov's theorem.* Given $k$ and $t$, and with $h(s) = \min(2s/t, 1)2kt^\alpha$ as defined above, we associate to a path $b$ a shifted path $b'$ by the relation

$$b'_s \equiv b_s - h(s) \qquad \text{for } s \in \mathbb{R}.$$

Notice that this shift transforms a path which lives in the interval $I_k$ for all $s \in [t/2, t]$ into a path which belongs to $I_0$ in the same time interval. More precisely, one immediately checks that $\mathbf{1}_{L_k}(b) = \mathbf{1}_{L_0}(b')$. Let us call $M_t(b')$ the Girsanov density involved in the shift between $b$ and $b'$, that is,

$$M_t(b') = \exp(-b'_{t/2}4kt^{\alpha-1} - 4k^2t^{2\alpha-1}).$$

The choice of $h(s) = 4kst^{\alpha-1}$ for $s \in [0, t/2]$ is made to obtain a continuous function that starts at 0, and is piecewise linear (constant over $[t/2, t]$); this function has the advantage that its Girsanov "energy" is minimal, ensuring that our proof is most efficient. It is possible that other, nonlinear, choices could have fulfilled our purposes, but this would be an unnecessary complication. For sake of clarity, let us stress now the dependence of the random variables $\delta, \eta$ and so on, on the data of our problem: it is readily checked, for instance, that

$$\eta_j = \eta_j(W) \quad \text{and} \quad \delta_j = \delta_j(b, \mathcal{L}(W)),$$

where a function of $(W)$ represents its dependence on the increments of $W$ in the interval $[0, t]$, as a random variable, where the symbol $\mathcal{L}(\cdot)$ denotes the law (distribution) of a process on $[0, t]$, and where a function of $b$ represents its dependence on the fixed path $b$. Then, adopting this convention, we have

$$\tilde{Z}_t^\alpha(k, W) = E_b\bigg[\mathbf{1}_{L_k} \exp\bigg(\beta \int_0^t W(ds, b'_s + h(s)) - \sum_{j \in \mathbb{Z}} \delta_j(b' + h, \mathcal{L}(W))\eta_j(W)\bigg)\bigg].$$

After applying Girsanov's transformation, noting that by definition, $\int_0^t W(ds, b'_s + h(s)) = \int_0^t W^{k,t}(ds, b'_s)$, we get (recall that $b'$ is a standard Brownian motion under the new probability, so that it is notationally legitimate to write $b$ instead of $b'$, and to denote expectation with respect to the new measure by $E_b$)

$$\tilde{Z}_t^\alpha(k,W) = E_b\left[\mathbf{1}_{L_0(b)} M_t(b) \exp\left(\beta\left(\int_0^t W^{k,t}(ds, b_s) r - \sum_{j\in\mathbb{Z}} \delta_j(b+h, \mathcal{L}(W))\eta_j(W)\right)\right)\right].$$

*Step* 2: *reexpressing the transformed $\eta$*. One should now compare the random variables $\eta_j(W)$ and $\eta_j(W^{k,t})$: by definition of these quantities, we have

$$\begin{aligned}
(6.4) \quad \eta_j(W^{k,t}) &= \frac{1}{t^{2\alpha}} \int_{(2j-1)t^\alpha}^{(2j+1)t^\alpha} \int_{t/2}^t W(ds, x+2kt^\alpha)\, dx \\
&= \frac{1}{t^{2\alpha}} \int_{(2(j+k)-1)t^\alpha}^{(2(j+k)+1)t^\alpha} \int_{t/2}^t W(ds, x)\, dx = \eta_{j+k}(W).
\end{aligned}$$

In particular, the law of $\eta(W^{k,t})$, considered as the set of random variables forming that sequence, is the same as the law of $\eta(W)$, a fact which we will not use in this proof, but will be crucial in the proof of the next lemma.

*Step* 3: *reexpressing the transformed $\delta$*. Along the same lines as (6.4), we now show that

(6.5) $$\delta_j(b+h, \mathcal{L}(W)) = \delta_{j-k}(b, \mathcal{L}(W^{k,t})).$$

To see this, we recall the definition of $\delta$: we have

$$\delta = \delta(b+h, \mathcal{L}(W)) = [C(t)]^{-1}v = [C(t, \mathcal{L}(W))]^{-1}v(b+h, \mathcal{L}(W)),$$

where we calculate

$$\begin{aligned}
&C_{\ell,m}(t, \mathcal{L}(W)) \\
&= \frac{1}{t^{(\alpha+1)}} \mathbf{E}\left[\int_{t/2}^t \int_{(2m-1)t^\alpha}^{(2m+1)t^\alpha} W(ds,x)\, dx \cdot \int_{t/2}^t \int_{(2\ell-1)t^\alpha}^{(2\ell+1)t^\alpha} W(ds,x)\, dx\right] \\
&= \frac{1}{t^{(\alpha+1)}} \mathbf{E}\left[\int_{t/2}^t \int_{(2(m-k)-1)t^\alpha}^{(2(m-k)+1)t^\alpha} W(ds, x+2kt^\alpha)\, dx \right.\\
&\qquad\qquad \left. \times \int_{t/2}^t \int_{(2(\ell-k)-1)t^\alpha}^{(2(\ell-k)+1)t^\alpha} W(ds, x+2kt^\alpha)\, dx\right] \\
&= \frac{1}{t^{(\alpha+1)}} \mathbf{E}\left[\int_{t/2}^t \int_{I_{m-k}} W^{k,t}(ds,x)\, dx \cdot \int_{t/2}^t \int_{I_{\ell-k}} W^{k,t}(ds,x)\, dx\right] \\
&= C_{\ell-k, m-k}(t, \mathcal{L}(W^{k,t})),
\end{aligned}$$

and similarly,

$$v_\ell(b + h, \mathcal{L}(W))$$
$$= 4t^{\alpha-1}\mathbf{E}\bigg[\int_{t/2}^t \int_{(2\ell-1)t^\alpha}^{(2\ell+1)t^\alpha} W(ds, x)\, dx \cdot \int_{t/2}^t W(ds, b_s + h(s))\bigg]$$
$$= 4t^{\alpha-1}\mathbf{E}\bigg[\int_{t/2}^t \int_{(2(\ell-k)-1)t^\alpha}^{(2(\ell-k)+1)t^\alpha} W(ds, x + h(s))\, dx \cdot \int_{t/2}^t W(ds, b_s + h(s))\bigg]$$
$$= 4t^{\alpha-1}\mathbf{E}\bigg[\int_{I_\ell} \int_{t/2}^t W^{k,t}(ds, x)\, dx \cdot \int_{t/2}^t W^{k,t}(ds, b_s)\bigg]$$
$$= v_{\ell-k}(b, \mathcal{L}(W^{k,t})).$$

We may thus write that the definition of $\delta(b + h, \mathcal{L}(W))$ is equivalent to

$$\forall \ell \in \mathbb{Z}\colon \sum_{m \in \mathbb{Z}} C_{\ell,m}(t, \mathcal{L}(W))\delta_m(b + h, \mathcal{L}(W)) = v_\ell(b + h, \mathcal{L}(W))$$
$$\iff \quad \forall \ell \in \mathbb{Z}\colon \sum_{m \in \mathbb{Z}} C_{\ell-k, m-k}(t, \mathcal{L}(W^{k,t}))\delta_m(b + h, \mathcal{L}(W))$$
$$= v_{\ell-k}(b, \mathcal{L}(W^{k,t}))$$
$$\iff \quad \forall \ell \in \mathbb{Z}\colon \sum_{m \in \mathbb{Z}} C_{\ell,m}(t, \mathcal{L}(W^{k,t}))\delta_{m+k}(b + h, \mathcal{L}(W))$$
$$= v_\ell(b, \mathcal{L}(W^{k,t})).$$

This last statement is equivalent to saying $\delta_{m+k}(b+h, \mathcal{L}(W)) = \delta_m(b, \mathcal{L}(W^{k,t}))$, which is precisely the statement of (6.5).

*Step* 4: *conclusion.* Plugging equations (6.4) and (6.5) into (6), we end up with

$$\tilde{Z}_t^\alpha(k, W) = E_b\bigg[\mathbf{1}_{L_0(b)} M_t(b) \exp\bigg(\beta\bigg(\int_0^t W^{k,t}(ds, b_s) \\ - \sum_{j \in \mathbb{Z}} \delta_{j-k}(b, \mathcal{L}(W^{k,t})) \eta_{j-k}(W^{k,t})\bigg)\bigg)\bigg].$$

To conclude the proof of the proposition, notice that for $b \in L_0$, we get $|b_{t/2}| \leq t^\alpha$, and therefore,

(6.6) $$M_t(b) \geq \exp(-4kt^{2\alpha-1} - 4k^2 t^{2\alpha-1}).$$

Combining (6) and (6.6), and renumbering the sum for $j \in \mathbb{Z}$ as $j' = j - k \in \mathbb{Z}$, we recognize the term $\tilde{Z}_t^\alpha(0, W^{k,t})$, and the proof is complete. $\square$

The above proof has an important consequence which we record here for use at a crucial point in the next section.

LEMMA 6.3. *Let*

$$X(W,b) = -H_t(b) - \sum_{j \in \mathbb{Z}} \delta_j \eta_j = \int_0^t W(ds, b_s) - \sum_{j \in \mathbb{Z}} \delta_j(b, \mathcal{L}(W)) \eta_j(W)$$

*and, therefore,*

$$X(W^{k,t}, b) = \int_0^t W^{k,t}(ds, b_s) - \sum_{j \in \mathbb{Z}} \delta_j(b, \mathcal{L}(W^{k,t})) \eta_j(W^{k,t}).$$

*Denote by $\eta(W)$ the entire sequence $\{\eta_j(W) : j \in \mathbb{Z}\}$. Then for each $b$, $X(W,b)$ and $\eta(W)$ are independent, and for each $k \in \mathbb{Z}$, and each $b$, $X(W^{k,t}, b)$ and $\eta(W)$ are independent.*

PROOF. We have already proved in Proposition 5.3(iv) that $X(W,b)$ and $\eta(W)$ are independent, which is the first half of what we have to prove. This implies in addition that $X(W^{k,t}, b)$ and $\eta(W^{k,t})$ are also independent because the random fields $W$ and $W^{k,t}$ have the same distribution (Lemma 6.1).

To conclude the proof this lemma, we simply invoke the portion of the proof of Proposition 6.2 which shows the specific shift equality relation $\eta_{j+k}(W) = \eta_j(W^{k,t})$, from (6.4): this is a **P**-almost-sure equality in $\Omega$. This implies that the sets of points in the sequences $\{\eta_j(W) : j \in \mathbb{Z}\}$ and $\{\eta_j(W^{k,t}) : j \in \mathbb{Z}\}$ are precisely the same sets of random variables. Therefore, for each $k$ and $b$, $X(W^{k,t}, b)$ is independent of the entire sequence $\eta(W)$. $\square$

**7. Proof of Lemma 3.3.** Recall that we have reduced our problem to the evaluation of $\mathbf{P}(\mathcal{B}_t)$, where

$$\mathcal{B}_t = \mathcal{A}_t^c = \{\text{For all } k \in \mathbb{Z}, Z_t^\alpha(k) \leq Z_t^\alpha(0)\},$$

and one wishes to show that $\lim_{t \to \infty} \mathbf{P}(\mathcal{B}_t) = 0$. Then a first step in order to prove this claim is to truncate $\mathcal{B}_t$: for a positive integer $M$, let $\mathbb{Z}_M$ and $\bar{\mathbb{Z}}_M$ be the sets defined respectively by

(7.1) $\quad \bar{\mathbb{Z}}_M = \{-M, -M+1, \ldots, M-1, M\} \quad \text{and} \quad \mathbb{Z}_M = \bar{\mathbb{Z}}_M \setminus \{0\},$

and $\mathcal{B}_{M,t}$ the event defined by

$$\mathcal{B}_{M,t} = \{\text{For all } k \in \mathbb{Z}_M, Z_t^\alpha(k) \leq Z_t^\alpha(0)\}.$$

Then obviously, $\mathbf{P}(\mathcal{B}_t) \leq \mathbf{P}(\mathcal{B}_{M,t})$, and we only need to prove that $\mathbf{P}(\mathcal{B}_{M,t})$ tends to 0 as $t \to \infty$.

Here is a brief account on the strategy we will follow in order to complete our proof.

(1) Recall that we are trying to bound

(7.2) $\mathbf{P}(\mathcal{B}_{M,t}) = \mathbf{P}(E_b[\mathbf{1}_{L_k} e^{-\beta H_t(b)}] < E_b[\mathbf{1}_{L_0} e^{-\beta H_t(b)}]$ for all $k \in \mathbb{Z}_M).$

A natural idea is then to split the conditions $E_b[\mathbf{1}_{L_k}e^{-\beta H_t(b)}] < E_b[\mathbf{1}_{L_0}e^{-\beta H_t(b)}]$ in terms of a condition involving the random variables $\eta_l$ introduced at (5.1), on which we have a reasonable control, and another set of conditions involving some random variables independent of the family $\{\eta_l; l \in \mathbb{Z}\}$. However, we have already seen in Proposition 5.3 that $-H_t(b) - \sum_{j \in \mathbb{Z}} \delta_j \eta_j$ is independent of $\{\eta_l; l \in \mathbb{Z}\}$. Thus, a natural choice will be to replace $e^{-\beta H_t(b)}$ by $e_t(b)$ in the expression (7.2), where $e_t(b)$ is defined by

$$e_t(b) := \exp\left(-\beta\left(H_t(b) + \sum_{j \in \mathbb{Z}} \delta_j \eta_j\right)\right).$$

Of course, this induces a correction term $\exp(\beta \sum_{j \in \mathbb{Z}} \delta_j \eta_j)$, but this term can be controlled, since the covariance structure of the family $\{\eta_l; l \in \mathbb{Z}\}$ is given by Proposition 4.2, and the vector $\delta$ is controlled by means of Proposition 5.3. Up to a negligible term, we will be allowed to bound $\mathbf{P}(\mathcal{B}_{M,t})$ by a probability of the form

(7.3) $\quad \mathbf{P}\left(\text{For any } k \in \mathbb{Z}_M; \frac{\tilde{Z}_t^\alpha(k)}{\tilde{Z}_t^\alpha(0)} < \exp(2\gamma t^{2\alpha-1} + \eta_k^*)\right),$

where $\tilde{Z}_t^\alpha(k) = E_b[\mathbf{1}_{L_k} e_t(b)]$, as was defined in Section 6 on Girsanov's theorem, the term $t^{2\alpha-1}$ comes from the sharp estimates of $\delta$ in Proposition 5.3, and the random variable $\eta_k^*$ is one which is defined using only the random variables $\eta$, because it results from using $e_t(b)$ instead of $e^{-H_t(b)}$. The effect of $\eta_k^*$ can be studied separately from the behavior of the ratio $\tilde{Z}_t^\alpha(k)/\tilde{Z}_t^\alpha(0)$, by the independence property of these two quantities.

(2) Notice that up to now, we have chosen our parameters carefully in order to get a penalization of order $\exp(2\gamma t^{2\alpha-1})$ in (7.3). This was chosen to be consistent with the correction $\exp(-4(k+k^2)t^{2\alpha-1})$ we must impose on $b$ if we wish that it live the second half of its life in $I_k$, as we showed by using Girsanov's theorem in Proposition 6.2. In fact, we will be able to bound $\mathbf{P}(\mathcal{B}_{M,t})$ by $\mathbf{P}(F_M)$, where the event $F_M$ is defined by

$$F_M = \left\{\text{For any } k \in \mathbb{Z}_M; \frac{\tilde{Z}_t^\alpha(0, W^{k,t})}{\tilde{Z}_t^\alpha(0, W)} < \exp(\hat{\gamma} t^{2\alpha-1} + \eta_k^*)\right\},$$

for some constant $\hat{\gamma} = \hat{\gamma}(M)$, where the shifted environments $W^{k,t}$ are defined in (6.1).

(3) It turns out that the random variable $\eta_k^*$ is optimally chosen to be of the order $\eta_0 - \eta_k$ [see the definition (7.9) we chose below]. We are now considering a set $F_M$ involving the random variables $\tilde{Z}_t$ and $\eta_k^*$, and this will allow us to take advantage of the following facts:

1. The ratio $\tilde{Z}_t^\alpha(0, W^{k,t})/\tilde{Z}_t^\alpha(0, W)$ cannot be too small at many different sites $k \in \mathbb{Z}_M$, by translation invariance in space of $W$.

2. Proposition 4.2 asserts that $\{t^{-(1-\alpha)/2}\eta_k; k \in \mathbb{Z}_M\}$ is asymptotically a standard Gaussian vector. Since $\eta_k^*$ is of the order $\eta_0 - \eta_k$ (and thus of magnitude $t^{(1-\alpha)/2}$), it can be highly negative at many different sites; thus, we are allowed to expect that $\exp(\hat{\gamma}t^{2\alpha-1} + \eta_k^*)$ is much smaller than 1 at many different sites of $\mathbb{Z}_M$.
3. The random variables $\tilde{Z}_t^\alpha$ are independent of anything defined using $\eta$, including $\eta_k^*$, and hence, the two effects alluded to above can be taken into account separately.

(4) These heuristic considerations will be formalized in step 3 of the proof below, through the introduction of an intricate family of subsets of $\mathbb{Z}_M$, but let us mention that the exponent $3/5$ comes out already at this stage: indeed, the above considerations only make sense if the magnitude $t^{(1-\alpha)/2}$ of the $\eta_k^*$ is greater than the magnitude $t^{2\alpha-1}$ of the penalization, so that a highly negative $\eta_k^*$ can win against the latter. This can only occur, obviously, whenever $\alpha < 3/5$. In this sense, our estimates are quite sharp: they mainly rely on the covariance structure of $\eta$ and on Girsanov's theorem applied to $b$.

Before going into the details of our calculations, let us introduce a new set $\mathcal{B}_{M,t}$: as mentioned above, our computations will bring out some expressions of the form $u_t := \sum_{j \in \mathbb{Z}} \delta_j \eta_j$, and it will be convenient to keep this kind of term of order $O(t^{2\alpha-1})$, which is also the order of the exponential correction term appearing in (6.3). However, since $\delta$ satisfies Proposition 5.3, it is easily checked that $u_t$ is of the desired order if $\eta_j \leq |j-k|^\tau t^{3\alpha-1}$ on $L_k$. These considerations motivate the introduction of the event

$$\mathcal{B}_{M,t} \equiv \{\text{There exists } \ell \in \bar{\mathbb{Z}}_M \text{ and } j \in \mathbb{Z}\setminus\{\ell\}; |\eta_j| \geq |j-\ell|^\tau t^{3\alpha-1}\},$$

and we will trivially bound $\mathbf{P}(\mathcal{B}_{M,t})$ by

$$(7.4) \qquad \mathbf{P}(\mathcal{B}_{M,t}) \leq \mathbf{P}(\mathcal{B}_{M,t}) + \mathbf{P}(\mathcal{B}_{M,t}^c \cap \mathcal{B}_{M,t}).$$

We will now prove that the two terms on the right-hand side of (7.4) vanish as $t \to \infty$, whenever $M$ is large enough.

*Step* 1: *estimation of* $\mathbf{P}(\mathcal{B}_{M,t})$. Let $\Phi$ be the distribution function of a standard Gaussian random variable, that is, if $Z \sim \mathcal{N}(0,1)$, then

$$(7.5) \qquad \Phi(x) = \mathbf{P}(Z \leq x),$$

and set $\bar{\Phi} = 1 - \Phi$. Then let us bound simply $\mathbf{P}(\mathcal{B}_{M,t})$ by

$$\mathbf{P}(\mathcal{B}_{M,t}) \leq \sum_{\ell \in \bar{\mathbb{Z}}_M} \sum_{j \neq \ell} \mathbf{P}(|\eta_j| \geq |j-\ell|^\tau t^{3\alpha-1})$$

$$\leq 2 \sum_{\ell \in \bar{\mathbb{Z}}_M} \sum_{j \neq \ell} \bar{\Phi}\left(\frac{2|j-\ell|^\tau t^{(7\alpha-3)/2}}{C_{0,0}^{1/2}(t)}\right),$$

where $C_{0,0}(t)$, defined in (4.2), equals $t^{\alpha-1}/4\mathbf{E}[\eta_\ell \eta_k]$. Recall that $\bar{\Phi}(x) \leq e^{-x^2/2}$ for $x$ large, enough, and that $C(t)$ satisfies Proposition 4.1. Thus, for two constants $c_1, c_2 > 0$, we get

$$\mathbf{P}(\mathcal{\tilde{B}}_{M,t}) \leq c_1 M \sum_{j \geq 1} \exp(-c_2 j^{2\tau} t^{7\alpha-3}). \tag{7.6}$$

The following facts are now easily seen:

- The series in the right-hand side of (7.6) is convergent, since $\tau > 0$, which explains the choice of the norm $\|x\|_{\tau,\ell}$ in order to bound $\eta_j$.
- Since we have assumed $\alpha > 1/2 > 3/7$, we have $7\alpha - 3 > 0$, and thus, an elementary application of the dominated convergence theorem yields

$$\lim_{t \to \infty} \mathbf{P}(\mathcal{\tilde{B}}_{M,t}) = 0, \tag{7.7}$$

which proves our first claim.

*Step* 2: *estimation of* $\mathbf{P}(\mathcal{\tilde{B}}^c_{M,t} \cap \mathcal{B}_{M,t})$. Recall that the vector $\delta$ has been introduced because $-H_t(b) - \sum_{j \in \mathbb{Z}} \delta_j \eta_j$ is independent of the family $\eta$, and for sake of compactness of notation, set

$$e_t(b) = \exp\left(-\beta\left(H_t(b) + \sum_{j \in \mathbb{Z}} \delta_j \eta_j\right)\right). \tag{7.8}$$

Now we have

$$\mathbf{P}(\mathcal{\tilde{B}}^c_{M,t} \cap \mathcal{B}_{M,t})$$
$$= \mathbf{P}(\mathcal{\tilde{B}}^c_{M,t} \text{ and } E_b[\mathbf{1}_{L_k} e^{-H_t(b)}] < E_b[\mathbf{1}_{L_0} e^{-H_t(b)}] \text{ for all } k \in \mathbb{Z}_M)$$
$$= \mathbf{P}\left(\mathcal{\tilde{B}}^c_{M,t} \text{ and } E_b\left[\mathbf{1}_{L_k} e_t(b) \exp\left(\sum_{j \in \mathbb{Z}} \beta \delta_j \eta_j\right)\right]\right.$$
$$\left. < E_b\left[\mathbf{1}_{L_0} e_t(b) \exp\left(\sum_{j \in \mathbb{Z}} \beta \delta_j \eta_j\right)\right] \text{ for all } k \in \mathbb{Z}_M\right).$$

As mentioned before, $\delta := C^{-1}(t) v$ depends on the path $b$, as is easily seen from definition (5.2). In order to get rid of the term $\sum_{j \in \mathbb{Z}} \delta_j \eta_j$, we will then set

$$\check{\eta}_0 = \max\left(\beta \underline{d} \eta_0, \beta \overline{d} \eta_0\right) \quad \text{and} \quad \hat{\eta}_k = \min\left(\beta \underline{d} \eta_k, \beta \overline{d} \eta_k\right), \tag{7.9}$$

where the constants $\underline{d}, \overline{d}$ have been introduced in Proposition 5.3. Then, according to the definition of $\mathcal{\tilde{B}}^c_{M,t}$, we get

$$\mathbf{P}(\mathcal{\tilde{B}}^c_{M,t} \cap \mathcal{B}_{M,t}) \leq \mathbf{P}\bigg(\text{For any } k \in \mathbb{Z}_M,$$

$$E_b\left[\mathbf{1}_{L_k}e_t(b)\exp\left(-\sum_{j\in\mathbb{Z}}\beta|\delta_j||j-k|^\tau t^{3\alpha-1}+\hat{\eta}_k\right)\right]$$

$$< E_b\left[\mathbf{1}_{L_0}e_t(b)\exp\left(\sum_{j\in\mathbb{Z}}\beta|\delta_j||j^\tau t^{3\alpha-1}+\check{\eta}_0\right)\right]).$$

Now, invoking Proposition 5.3, item (iii), we obtain that, for any integer $k$, there exists a constant $\gamma$ (possibly depending on $\beta$) such that $\sum_{j\in\mathbb{Z}}\beta|\delta_j||j-\ell|^\tau \leq \gamma t^{-\alpha}$ on $L_k$. Thus, thanks to the fact that the random variables $\eta$ only depend on $W$, and observing that $\tilde{Z}_t^\alpha(k) = E_b[\mathbf{1}_{L_k}e_t(b)]$, we get

(7.10)
$$\mathbf{P}(\mathcal{B}_{M,t}^c \cap \mathcal{B}_{M,t})$$
$$\leq \mathbf{P}(\text{For any } k \in \mathbb{Z}_M;$$
$$\tilde{Z}_t^\alpha(k)\exp(-\gamma t^{2\alpha-1}+\hat{\eta}_k) < \exp(\gamma t^{2\alpha-1}+\check{\eta}_0)\tilde{Z}_t^\alpha(0))$$
$$= \mathbf{P}\left(\text{For any } k \in \mathbb{Z}_M; \frac{\tilde{Z}_t^\alpha(k)}{\tilde{Z}_t^\alpha(0)} < \exp(2\gamma t^{2\alpha-1}+\check{\eta}_0 - \hat{\eta}_k)\right).$$

Let us apply now Proposition 6.2 in order to conclude that

$$\mathbf{P}(\mathcal{B}_{M,t}^c \cap \mathcal{B}_{M,t})$$
$$\leq \mathbf{P}\left(\text{For any } k \in \mathbb{Z}_M; \frac{\tilde{Z}_t^\alpha(0,W^{k,t})}{\tilde{Z}_t^\alpha(0,W)} < \exp(\hat{\gamma}t^{2\alpha-1}+\check{\eta}_0 - \hat{\eta}_k)\right),$$

where $\hat{\gamma} = \hat{\gamma}(M) = \sup\{2\gamma + \zeta(k); k \in \mathbb{Z}_M\}$ and $\zeta(k) = 4k(k+1)$. We have thus proved that

$$\mathbf{P}(\mathcal{B}_{M,t}^c \cap \mathcal{B}_{M,t}) \leq \mathbf{P}(F_M),$$

where

$$F_M = \left\{\text{For any } k \in \mathbb{Z}_M; \frac{\tilde{Z}_t^\alpha(0,W^{k,t})}{\tilde{Z}_t^\alpha(0,W)} < \exp(\hat{\gamma}t^{2\alpha-1}+\check{\eta}_0 - \hat{\eta}_k)\right\}.$$

*Step* 3: *evaluation of* $\mathbf{P}(F_M)$. We can see now that the probability of $F_M$ will be expressed in terms of a balance between the values of $\check{\eta}_0 - \hat{\eta}_k$ (which will be assumed to be highly negative) and the ratio $\tilde{Z}_t^\alpha(0,W^{k,t})/\tilde{Z}_t^\alpha(0,W)$, which cannot be too small at many different sites $k$. In order to quantify this heuristic statement, we introduce a family $\vec{\mathcal{S}}_{M,m}$ of subsets of $\bar{\mathbb{Z}}_M$ which will be used to construct a large symmetric set $L$ around 0 such that $\check{\eta}_0 - \hat{\eta}_\ell < -t^{2\alpha-1+\rho}$ for all $\ell \in L$: for a given $\rho > 0$ and integer numbers $m$ and $M$, define the families of subsets

$$\mathcal{S}_{M,m} = \bigcup_{k,\hat{k}\in D_{M,m}}\{k\mathbb{Z}_{\hat{k}}\},$$

(7.11) $$\text{with } D_{M,m} = \{(k,k'): k \geq 1, \hat{k} \geq m; k\mathbb{Z}_{\hat{k}} \subset \mathbb{Z}_M\}$$

$$\breve{\mathcal{S}}_{M,m} = \{L \subset \mathbb{Z}_M; \text{There exists } S \in \mathcal{S}_{M,m} \text{ such that } S \subset L\}.$$

In relation with these families of subsets of $\mathbb{Z}_M$, set also

(7.12) $$\hat{F}_{M,m,\rho} = \bigcup_{L \in \breve{\mathcal{S}}_{M,m}} \hat{F}_{\rho,L},$$

with

(7.13) $$\hat{F}_{\rho,L} = \{\breve{\eta}_0 - \hat{\eta}_\ell < -t^{2\alpha-1+\rho}, \text{for all } \ell \in L,$$
$$\breve{\eta}_0 - \hat{\eta}_{\hat{\ell}} > -t^{2\alpha-1+\rho}, \text{for all } \hat{\ell} \in \mathbb{Z}_M \setminus L\}.$$

Then one can bound trivially $\mathbf{P}(F_M)$ by

$$\mathbf{P}(F_M) \leq 1 - \mathbf{P}(\hat{F}_{M,m,\rho}) + \mathbf{P}(F_M \cap \hat{F}_{M,m,\rho}).$$

Furthermore, for $t$ large enough, we have $\hat{\gamma} t^{2\alpha-1} - t^{2\alpha-1+\rho} < 0$, which explains the need for the constant $\rho > 0$. Thus,

$$F_M \cap \hat{F}_{M,m,\rho} \subseteq \bigcup_{L \in \breve{\mathcal{S}}_{M,m}} \bigcap_{\ell \in L} \left\{ \frac{\tilde{Z}^\alpha_t(0, W^{\ell,t})}{\tilde{Z}^\alpha_t(0,W)} < \exp(\hat{\gamma} t^{2\alpha-1} - 2t^{2\alpha-1+\rho}) \right\} \cap \hat{F}_{\rho,L}$$

$$\subseteq \bigcup_{L \in \breve{\mathcal{S}}_{M,m}} \{\tilde{Z}^\alpha_t(0, W^{\ell,t}) < \tilde{Z}^\alpha_t(0,W) \text{ for all } \ell \in L\} \cap \hat{F}_{\rho,L}.$$

Hence, we get

(7.14) $$\mathbf{P}(F_M) \leq 1 - \mathbf{P}(\hat{F}_{M,m,\rho})$$
$$+ \sum_{L \in \breve{\mathcal{S}}_{M,m}} \mathbf{P}(\{\tilde{Z}^\alpha_t(0, W^{\ell,t}) < \tilde{Z}^\alpha_t(0,W) \text{ for all } \ell \in L\} \cap \hat{F}_{\rho,L})$$
$$\leq 1 - \mathbf{P}(\hat{F}_{M,m,\rho})$$
$$+ \sum_{L \in \breve{\mathcal{S}}_{M,m}} \mathbf{P}(\tilde{Z}^\alpha_t(0, W^{\ell,t}) < \tilde{Z}^\alpha_t(0,W) \text{ for all } \ell \in L) \mathbf{P}(\hat{F}_{\rho,L}),$$

where in the last step, we have used the independence, proved in the next step, between the random variables $\tilde{Z}^\alpha_t(0, W^{\ell,t})$ and the sequence $\{\eta_k; k \in \bar{\mathbb{Z}}_M\}$, and also between $\tilde{Z}^\alpha_t(0,W)$ and the sequence $\{\eta_k; k \in \bar{\mathbb{Z}}_M\}$.

*Step 4*: *independence of $\eta$ and the $\tilde{Z}^\alpha_t$'s*. Using the notation $X(W,b)$ introduced in Lemma 6.3, this lemma's conclusion is that $X(W,b)$ and $\eta(W)$ are independent for each continuous function $b$; after evaluation of $\tilde{Z}^\alpha_t(0,W)$ in formula (6.2), it implies that the latter is also independent of $\eta$.

Lemma 6.3 can also be applied to prove the other independence: it proves that for each fixed $b, k$, we have independence of $X(W^{k,t}, b)$ and the entire sequence $\eta$. When defining $\tilde{Z}_t^\alpha(0, W^{\ell,t})$, formula (6.2) must be used with $W$ replaced by $W^{\ell,t}$, which specifically means

$$\tilde{Z}_t^\alpha(0, W^{\ell,t}) = E_b\left[\mathbf{1}_{L_k} \exp\beta\left(\int_0^t W^{\ell,t}(ds, b_s) - \sum_{j \in \mathbb{Z}} \delta_j(b, \mathcal{L}(W^{\ell,t}))\eta_j(W^{\ell,t})\right)\right]$$

$$= E_b[\mathbf{1}_{L_k} \exp\beta(X(W^{\ell,t}, b))],$$

proving that $\tilde{Z}_t^\alpha(0, W^{\ell,t})$ is independent of $\eta$, as required to justify (7.14) in step 3.

One can prove in addition that $\delta(b, \mathcal{L}(W)) = \delta(b, \mathcal{L}(W^{\ell,t}))$ for any $\ell$, but this fact will not be needed.

*Step* 5: *finishing the proof.* The end of our proof of Lemma 3.3 relies on the following propositions, whose proofs will be postponed until the next sections.

PROPOSITION 7.1. *Let $m$ be a fixed positive even integer, and $M > m$. Then, for any $L \in \breve{\mathscr{S}}_{M,m}$, we have*

$$\mathbf{P}(\tilde{Z}_t^\alpha(0, W^{\ell,t}) < \tilde{Z}_t^\alpha(0, W) \text{ for all } \ell \in L) \leq \frac{1}{m}.$$

PROPOSITION 7.2. *Let $m$ be a fixed positive integer. Let $\rho$ be a strictly positive number such that $\frac{5}{2}(\alpha - \frac{3}{5}) + \rho < 0$. Then, for $t$ large enough, there exists a $M$ large enough such that*

(7.15) $$\mathbf{P}(\hat{F}_{M,m,\rho}) \geq 1 - \frac{1}{m}.$$

With these results in mind, let us finish now the proof of Lemma 3.3, and thus of our theorem: take $t, M$ large enough so that (7.15) is satisfied. Then (7.14) yields directly, invoking Proposition 7.1 and the fact that the events $\hat{F}_{\rho,L}$ are disjoints,

$$\mathbf{P}(F_M) \leq \frac{1}{m} + \frac{1}{m} \sum_{L \in \breve{\mathscr{S}}_{M,m}} \mathbf{P}(\hat{F}_{\rho,L}) \leq \frac{1}{m} + \frac{1}{m} = \frac{2}{m},$$

which tends to 0 as $m \to \infty$, and ends the proof of the theorem, modulo establishing the last two propositions above.

Before proceeding with the proofs of Propositions 7.1 and 7.2, we discuss the consequences of weakening Hypothesis 3.1. If we assume only that

(7.16) $$Q(x) \leq |x|^{-2-\theta},$$

can we find values of $\theta \leq 1$ such that we still get superdiffusive behavior for the polymer, that is, $\alpha > 1/2$? Since the result of the Girsanov theorem,

Proposition 6.2, is not effected by the value of $\theta$ above, this means that the penalization from Girsanov's theorem, of order $t^{2\alpha-1}$, cannot be made smaller by a different choice of decorrelation speed in $Q$. Therefore, we should expect not to be able to preserve the threshold $\alpha < 3/5$. To see exactly what happens to this threshold under condition (7.16), we first state, and leave it to the reader to check, that we can rework the proof of Proposition 5.3, item (iii), to obtain instead

$$|\delta|_{\tau,k} - \delta_k = o(t^{-\alpha\theta}).$$

It is then simple to check that (7.10) becomes

$$\mathbf{P}(\mathcal{B}^c_{M,t} \cap \mathcal{B}_{M,t}) \leq \mathbf{P}\bigg(\text{For any } k \in \mathbb{Z}_M; \frac{\tilde{Z}^\alpha_t(k)}{\tilde{Z}^\alpha_t(0)} < \exp(2\gamma t^{3\alpha-1-\theta} + \check{\eta}_0 - \hat{\eta}_k)\bigg).$$

Hence, the application of Proposition 6.2 still works, but we can no longer make the corresponding Girsanov penalization of the same order, since for $\theta < 1$, $3\alpha - 1 - \alpha\theta > 2\alpha - 1$. Having thus convinced ourselves that Hypothesis 3.1 is the only way to get the entire proof to be efficient in terms of using comparable penalizations throughout, we can now ignore this inefficiency, and answer the question at the beginning of this paragraph. The reader will check that any other occurrences of the use of Hypothesis 3.1 are not further effected by switching to (7.16): the entire proof can still be used if we only require that the magnitude of the $\eta_k$'s, namely, $t^{(1-\alpha)/2}$, is larger than the new penalization $t^{3\alpha-1-\alpha\theta}$. This yields

$$\alpha < \frac{3}{7 - 2\theta}.$$

Now we see that to get a super-diffusive behavior, we need $3/(7-2\theta) > 1/2$, that is, $\theta > 1/2$. We also see that the weakest hypothesis required for such behavior is $Q(x) \leq x^{-5/2-\vartheta}$ for $\vartheta > 0$. We state these findings formally, using the reparametrization $\theta = \vartheta + 1/2$.

COROLLARY 7.3.  *Assume instead of Hypothesis 3.1 that there exists $\vartheta \in (0, 1/2]$ such that as $|x| \to \infty$,*

$$Q(x) = O(|x|^{-5/2-\vartheta}).$$

*Then for any $\varepsilon > 0$, we obtain the following specific super-diffusive behavior for the polymer measure:*

$$\lim_{t \to \infty} \mathbf{P}\bigg[\bigg\langle \sup_{s \leq t} |b_s| \bigg\rangle_t \geq t^{1/2 + \vartheta/(6-2\vartheta) - \varepsilon}\bigg] = 1.$$

7.1. *Proof of Proposition 7.1.* Let $L \in \breve{\mathscr{S}}_{M,m}$. Then, by definition (7.11) of $\breve{\mathscr{S}}_{M,m}$, there exists $k \geq 1$ such that $k\mathbb{Z}_m \subset L$. Then

$$\mathbf{P}(\tilde{Z}_t^\alpha(0, W^{\ell,t}) < \tilde{Z}_t^\alpha(0, W) \text{ for all } \ell \in L)$$
$$\leq \mathbf{P}(\tilde{Z}_t^\alpha(0, W^{\ell,t}) < \tilde{Z}_t^\alpha(0, W) \text{ for all } \ell \in k\mathbb{Z}_m).$$

It is thus sufficient to estimate the right-hand side in the above inequality.

Given an even integer $m \leq M$, recall that $\bar{\mathbb{Z}}_m$ has been defined at (7.1). Set also $\hat{m} = m/2$, and for each $i \in k\bar{\mathbb{Z}}_{\hat{m}}$, we associate the following event:

$$\Omega^{(i)} \equiv \{\tilde{Z}_t^\alpha(0, W^{\ell,t}) < \tilde{Z}_t^\alpha(0, W^{i,t}) \text{ for all } \ell \in k\bar{\mathbb{Z}}_{\hat{m}} \setminus \{i\}\}.$$

Then these events are disjoint, and since $|k\bar{\mathbb{Z}}_{\hat{m}}| = 2\hat{m} + 1$, we get trivially the existence of $i_0 \in k\bar{\mathbb{Z}}_{\hat{m}}$ such that

$$(7.17) \qquad \mathbf{P}(\Omega^{(i_0)}) \leq \frac{1}{2\hat{m}+1} \leq \frac{1}{m}.$$

However, the translation-invariance of the environment $W$ yields

$$\mathbf{P}(\tilde{Z}_t^\alpha(0, W^{\ell,t}) < \tilde{Z}_t^\alpha(0, W) \text{ for all } \ell \in k\mathbb{Z}_m)$$
$$= \mathbf{P}(\tilde{Z}_t^\alpha(0, W^{\ell+i_0,t}) < \tilde{Z}_t^\alpha(0, W^{i_0,t}) \text{ for all } \ell \in k\mathbb{Z}_m).$$

Indeed, exactly as we proved Lemma 6.1, denoting again $h(s) = \min(2s/t, 1)2kt^\alpha$, it holds that, for fixed $b$, $\int_0^t W(ds, b_s + (\ell + i_0)h(s))$ has the same distribution as $\int_0^t W(ds, b_s + \ell h(s))$.

We may now rewrite the above expression as the following upper bound:

$$(7.18) \quad \begin{aligned} &\mathbf{P}(\tilde{Z}_t^\alpha(0, W^{\ell,t}) < \tilde{Z}_t^\alpha(0, W) \text{ for all } \ell \in k\mathbb{Z}_m) \\ &\leq \mathbf{P}(\tilde{Z}_t^\alpha(0, W^{\ell,t}) < \tilde{Z}_t^\alpha(0, W^{i_0,t}) \text{ for all } \ell \in k\bar{\mathbb{Z}}_{\hat{m}} \setminus \{i_0\}) = \mathbf{P}(\Omega^{(i_0)}). \end{aligned}$$

Observe that the last inequality is just due to the elementary fact that $k\bar{\mathbb{Z}}_{\hat{m}} \setminus \{i_0\} \subset i_0 + k\mathbb{Z}_m$ whenever $i_0 \in k\bar{\mathbb{Z}}_{\hat{m}}$, a fact which is easily checked. Hence, putting together (7.17) and (7.18), we get the announced result.

7.2. *Proof of Proposition 7.2.* Recall that $\hat{F}_{M,m,\rho}$ is defined by (7.12), and define the quantity

$$\tau(t) := 2\beta^{-1} t^{(5/2)(\alpha - 3/5) + \rho},$$

which tends to 0 as $t \to \infty$ if $\alpha < \frac{3}{5}$ and $\rho$ is small enough. The following inequality

$$(7.19) \quad \mathbf{P}(\hat{F}_{M,m,\rho}) \geq \mathbf{P}\left(\bigcup_{L \in \breve{\mathscr{S}}_{M,m}} \{t^{(\alpha-1)/2}(\check{\eta}_0 - \hat{\eta}_\ell) \leq -\beta \tau(t) \text{ for all } \ell \in L\}\right)$$

is then easily established by an elementary inclusion argument, which we detail here. Indeed, assume that, for some $L \in \acute{\mathscr{S}}_{M,m}$, for all $\ell \in L$, $\eta$ satisfies

$$t^{(\alpha-1)/2}(\check{\eta}_0 - \hat{\eta}_\ell) \leq -\beta\tau(t),$$

which is equivalent to

$$\check{\eta}_0 - \hat{\eta}_\ell \leq -t^{2\alpha-1+\rho}.$$

To justify the above inequality, we only need to prove that, for some other $L' \in \acute{\mathscr{S}}_{M,m}$, the same $\eta$ also satisfies the above inequality for all $\ell \in L'$, while for all $\ell \in \mathbb{Z}_M \setminus L'$, the contrary holds, namely,

$$\check{\eta}_0 - \hat{\eta}_\ell > -t^{2\alpha-1+\rho}.$$

Let then $\Lambda$ be the subset of $\mathbb{Z}_M$ defined by

$$\Lambda = \{\ell \in \mathbb{Z}_M; \check{\eta}_0 - \hat{\eta}_\ell > -t^{2\alpha-1+\rho}\},$$

and set $L' = \mathbb{Z}_M \setminus \Lambda$. Then, by construction, $L'$ has the required properties defined above, and since $L' \supset L$, by definition of $\acute{\mathscr{S}}_{M,m}$, we have $L' \in \acute{\mathscr{S}}_{M,m}$.

In order to get a lower bound on the right-hand side of (7.19) above, we will construct now a large enough collection of symmetric and disjoint sets in $\mathbb{Z}_M$: with $m < M$, consider the collection $\{Q_q(m)\mathbb{Z}_m; q < q^*\}$, where the integers $Q_q(m)$ are defined by

$$Q_1(m) = 1, \qquad Q_{q+1}(m) = mQ_q(m) + 1, \qquad q^* = \inf\{q; Q_q(m) > M\}.$$

This collection is the sequence

$$\mathbb{Z}_m, (m+1)\mathbb{Z}_m, [m(m+1)+1]\mathbb{Z}_m, \ldots, Q_q(m)\mathbb{Z}_m, \ldots, Q_{q^*-1}(m)\mathbb{Z}_m,$$

which are nonoverlapping annuli in $\mathbb{Z}_M$, and therefore are indeed symmetric and disjoint subsets of $\mathbb{Z}_M$. Since $Q_q(m)\mathbb{Z}_m$ is certainly of the form $k\mathbb{Z}_{\hat{k}}$ with $k \geq 1$ and $\hat{k} \geq m$, and is a subset of $\mathbb{Z}_M$ as soon as $q < q^*$, by definition, $Q_q(m)\mathbb{Z}_m \in \acute{\mathscr{S}}_{M,m}$. Thus, using the notation $\check{\eta}_0, \hat{\eta}_\ell$ and $\eta_\ell$ defined in (7.8) and (7.9), and reverting to the notation $\tilde{\eta} = 2t^{-(1-\alpha)/2}\eta$, we get

$$\mathbf{P}(\hat{F}_{M,m,\rho}) \geq \mathbf{P}\left(\bigcup_{q<q^*} \{\max(\underline{d}\tilde{\eta}_0, \overline{d}\tilde{\eta}_0) - \min(\underline{d}\tilde{\eta}_\ell, \overline{d}\tilde{\eta}_\ell) \leq -\tau(t)$$
$$\text{for all } \ell \in Q_q(m)\mathbb{Z}_m\}\right).$$

Indeed, the original set $\hat{F}_{M,m,\rho}$ defined in (7.12) and (7.13) was a union of events indexed by $L \in \acute{\mathscr{S}}_{M,m}$, while here we use only sets of the form $L = Q_q(m)\mathbb{Z}_m$; moreover, the above condition on the difference $\max(\underline{d}\tilde{\eta}_0, \overline{d}\tilde{\eta}_0) - \min(\underline{d}\tilde{\eta}_\ell, \overline{d}\tilde{\eta}_\ell)$ is implied by the two conditions on the individual terms of this difference in $\hat{F}_{M,m,\rho}$, and the shorthand notation $\tau(t)$ was introduced above to be consistent with these conditions in (7.13). Let us call now $A_\ell$ the event

$$A_\ell = \{\max(\underline{d}\tilde{\eta}_0, \overline{d}\tilde{\eta}_0) - \min(\underline{d}\tilde{\eta}_\ell, \overline{d}\tilde{\eta}_\ell) \leq -\tau(t)\},$$

and we distinguish two cases according to the values of $\tilde{\eta}_0$:

(a) If $\tilde{\eta}_0 \geq 0$, then $\max(\underline{d}\tilde{\eta}_0, \overline{d}\tilde{\eta}_0) = \overline{d}\tilde{\eta}_0$, and hence, $A_\ell$ is the event defined by the relation

$$\min(\underline{d}\tilde{\eta}_\ell, \overline{d}\tilde{\eta}_\ell) \geq \tau(t) + \overline{d}\tilde{\eta}_0.$$

In particular, $\tilde{\eta}_\ell$ has to be positive, and thus, $A_\ell$ can be written as

$$\{\overline{d}\tilde{\eta}_0 - \underline{d}\tilde{\eta}_\ell < -\tau(t)\}.$$

(b) If $\tilde{\eta}_0 \leq -\tau(t)/\underline{d} \leq 0$, then $\max(\underline{d}\tilde{\eta}_0, \overline{d}\tilde{\eta}_0) = \underline{d}\tilde{\eta}_0$. Thus, $A_\ell$ can be written as the event defined by the relation

$$\min(\underline{d}\tilde{\eta}_\ell, \overline{d}\tilde{\eta}_\ell) \geq \tau(t) + \underline{d}\tilde{\eta}_0, \tag{7.20}$$

and if $\tilde{\eta}_0 \leq -\tau(t)/\underline{d}$, we have $\tau(t) + \underline{d}\tilde{\eta}_0 \leq 0$. Hence, (7.20) is implied by $\tilde{\eta}_\ell \geq 0$.

Summarizing the considerations above, we get

$$\mathbf{P}(\hat{F}_{M,m,\rho}) \geq \mathbf{P}(D^+) + \mathbf{P}(D^-),$$

with

$$D^+ = \bigcup_{q<q^*} \{\overline{d}\tilde{\eta}_0 - \underline{d}\tilde{\eta}_\ell \leq -\tau(t) \text{ for all } \ell \in Q_q(m)\mathbb{Z}_m\} \cap \{\tilde{\eta}_0 > 0\},$$

$$D^- = \bigcup_{q<q^*} \{\tilde{\eta}_\ell \geq 0 \text{ for all } \ell \in Q_q(m)\mathbb{Z}_m\} \cap \{\tilde{\eta}_0 \leq -\tau(t)/\underline{d}\}.$$

We will now prove that $\mathbf{P}(D^+)$ is close to $1/2$. Entirely similar arguments, left to the reader, lead to showing that $\mathbf{P}(D^-)$ can also be made arbitrarily close to $1/2$, concluding the proof of the proposition.

Observe that, according to Proposition 4.1, the random variables $\{\tilde{\eta}_\ell; l \in \bar{\mathbb{Z}}_M\}$ converge in distribution to a family of independent standard Gaussian random variables $\{\Upsilon_\ell; l \in \bar{\mathbb{Z}}_M\}$. Consequently, and using the fact that $-\tau(t) \to 0$ as $t \to \infty$,

$$\mathbf{P}(D^+) = \mathbf{P}\left(\bigcup_{q<q^*}\{\overline{d}\Upsilon_0 - \underline{d}\Upsilon_\ell \leq 0 \text{ for all } \ell \in Q_q(m)\mathbb{Z}_m\} \cap \{\Upsilon_0 > 0\}\right) + \varepsilon_M(t),$$

where, for a fixed $M \in \mathbb{N}$, we have $\lim_{t \to \infty} \varepsilon_M(t) = 0$. Furthermore, since the $\Upsilon_\ell$ are independent random variables, we get

$$\mathbf{P}(D^+) = \int_0^\infty \mathbf{P}\left(\bigcup_{q < q^*} \{\overline{d}x - \underline{d}\Upsilon_\ell \leq 0 \text{ for all } \ell \in Q_q(m)\mathbb{Z}_m\}\right)$$

(7.21)
$$\times \frac{e^{-x^2/2}}{(2\pi)^{1/2}} \, dx + \varepsilon_M(t)$$

$$= \frac{1}{2} - \int_0^\infty \mathbf{P}\left(\bigcap_{q < q^*} \hat{D}_q\right) \frac{e^{-x^2/2}}{(2\pi)^{1/2}} \, dx + \varepsilon_M(t),$$

where

$$\hat{D}_q = \{\text{There exists } \ell \in Q_q(m)\mathbb{Z}_m; \overline{d}x - \underline{d}\Upsilon_\ell \geq 0\}.$$

In order to take advantage of the independence of the $\Upsilon_\ell$, it is convenient to pick some disjoint sets out of $\mathbb{Z}_M$, which explains the choice of disjoint subsets $Q_q(m)\mathbb{Z}_m$. Now, it is easily seen that, for a fixed value $q_0$, if one desires to have $q^* > q_0$, it is sufficient to take $M$ of order $m^{q_0}$. Let us assume that we are in this situation; this means that, setting $\kappa = \overline{d}/\underline{d}$, we have

$$\mathbf{P}\left(\bigcap_{q < q^*} \hat{D}_q\right) \leq \mathbf{P}\left(\bigcap_{q \leq q_0} \{\text{There exists } \ell \in Q_q(m)\mathbb{Z}_m; \Upsilon_\ell \leq \kappa x\}\right)$$

$$= \mathbf{P}^{q_0}(\text{There exists } \ell \in \mathbb{Z}_m; \Upsilon_\ell \leq \kappa x) = [1 - \mathbf{P}^{2m}(\Upsilon_1 \geq \kappa x)]^{q_0}.$$

Plugging these inequalities into (7.21), we obtain

$$\mathbf{P}(D^+) \geq \frac{1}{2} - \int_0^\infty [1 - \mathbf{P}^{2m}(\Upsilon_1 \geq \kappa x)]^{q_0} \frac{e^{-x^2/2}}{(2\pi)^{1/2}} \, dx + \varepsilon_M(t).$$

Recall that the function $\Phi$ has been defined by relation (7.5). Then the last inequality yields

$$\mathbf{P}(D^+) \geq \frac{1}{2} - \int_0^\infty [1 - \Phi(\kappa x)^{2m}]^{q_0} \frac{e^{-x^2/2}}{(2\pi)^{1/2}} \, dx + \varepsilon_M(t).$$

It is easily seen that this probability can be made as close as we wish to $\frac{1}{2}$ by taking $q_0 \to \infty$, because $1/2 \leq \Phi(x) < 1$ for all $x \geq 0$, this asymptotic being equivalent to $M \to \infty$.

S. Bezerra
S. Tindel
Institut Elie Cartan
Université de Nancy 1 BP 239
54506-Vandoeuvre-lès-Nancy
France
E-mail: bezerra@iecn.u-nancy.fr
tindel@iecn.u-nancy.fr

F. Viens
Departments of Statistics and Mathematics
Purdue University
150 N. University St.
West Lafayette, Indiana 47907-2067
USA
E-mail: viens@purdue.edu